\documentclass[a4paper,oneside]{amsart}

\usepackage[a4paper]{geometry}
\usepackage[font=scriptsize]{caption}
\usepackage{bbm}
\usepackage{graphicx}
\usepackage{amssymb}
\usepackage{wrapfig}
\usepackage{float}
\usepackage{array}

\usepackage{xcolor}
\usepackage{etoolbox}
\usepackage[absolute]{textpos}

\usepackage[]{hyperref}
\definecolor{hyperrefLinkColor}{rgb}{0.1,0.1,1.0}
\definecolor{hyperrefCiteColor}{rgb}{0.75,0,0}
\definecolor{hyperrefUrlColor}{rgb}{0,0,0.8}
\definecolor{hyperrefFileColor}{rgb}{0.8,0.1,0.8}
\hypersetup{
    pdftitle={CMC surfaces with genus},
    pdfauthor={A. Bobenko, S. Heller, N. Schmitt},
    colorlinks=true,
    urlcolor=hyperrefUrlColor,
    filecolor=hyperrefFileColor,
    linkcolor=hyperrefLinkColor,
    citecolor=hyperrefCiteColor
}
\usepackage{cleveref}

\numberwithin{equation}{section}

\crefname{figure}{figure}{figures}

\newcommand{\MatrixGroup}[1]{{\mathrm{#1}}}

\newcommand{\matGL}[2]{\MatrixGroup{GL}_{#1}{#2}}

\newcommand{\matSL}[2]{\MatrixGroup{SL}_{#1}{#2}}
\newcommand{\matSO}[2]{\MatrixGroup{SO}_{#1}{#2}}
\newcommand{\matSU}[2]{\MatrixGroup{SU}_{#1}{#2}}
\newcommand{\matSpin}[2]{\MatrixGroup{Spin}_{#1}{#2}}

\newcommand{\matgl}[2]{\mathfrak{gl}({#1},\,{#2})}
\newcommand{\matsl}[2]{\mathfrak{sl}({#1},\,{#2})}

\newcommand{\matsu}[2]{\mathfrak{su}_{#1}{#2}}
\newcommand{\matsym}[2]{\MatrixGroup{Sym}_{#1}{#2}}
\newcommand{\matN}[2]{\MatrixGroup{M}_{#1\times #1}{#2}}

\makeatletter

\def\part{\@startsection{part}{0}%
  \z@{\linespacing\@plus\linespacing}{.5\linespacing}%
  {\normalfont\Large\bfseries\raggedright}}

\def\section{\@startsection{section}{1}%
  \z@{.7\linespacing\@plus\linespacing}{.5\linespacing}%
  {\normalfont\large\bfseries\centering}}

\makeatother

\newcommand{\transpose}[1]{{{#1}^{\top}}}

\newcommand{\half}{\tfrac{1}{2}}

\newcommand{\fourth}{\tfrac{1}{4}}

\newcommand{\ol}[1]{\overline{#1}}

\DeclareMathSymbol{\varnothing}{\mathord}{AMSb}{"3F}

\DeclareMathOperator{\adj}{adj}

\DeclareMathOperator{\del}{\partial\!}
\DeclareMathOperator{\diag}{diag}

\DeclareMathOperator{\id}{\mathbbm{1}}
\DeclareMathOperator{\image}{image}

\DeclareMathOperator{\rank}{rank}
\DeclareMathOperator*{\res}{res}

\DeclareMathOperator{\suchthat}{|}

\DeclareMathOperator{\tr}{tr}
\DeclareMathOperator{\trace}{tr}

\newcommand{\bbC}{\mathbb{C}}

\newcommand{\bbH}{\mathbb{H}}

\newcommand{\bbN}{\mathbb{N}}

\newcommand{\bbR}{\mathbb{R}}
\newcommand{\bbS}{\mathbb{S}}

\newcommand{\bbZ}{\mathbb{Z}}

\newcommand{\bbi}{\mathbbm{i}}

\newcommand{\calD}{\mathcal{D}}

\newcommand{\calG}{\mathcal{G}}

\newcommand{\calL}{\mathcal{L}}

\newcommand{\calP}{\mathcal{P}}

\newcommand{\CPone}{\bbC\mathrm{P}^1}

\newcommand{\deriv}{\mathrm{d}}

\newcommand{\Rstar}{{\bbR^\times}}

\newcommand*{\coloneq}{
  \mathrel{\vcenter{\baselineskip0.5ex \lineskiplimit0pt
      \hbox{\scriptsize.}\hbox{\scriptsize.}}}%
  =}

\newcommand{\spacecomma}{\,\,\,,}
\newcommand{\spaceperiod}{\,\,\,.}

\theoremstyle{plain}

\newtheorem{theorem}{Theorem}[section]
\newtheorem*{theorem*}{Theorem}

\newtheorem{lemma}[theorem]{Lemma}
\newtheorem*{lemma*}{Lemma}
\newtheorem{corollary}[theorem]{Corollary}
\newtheorem*{corollary*}{Corollary}

\newtheorem{proposition}[theorem]{Proposition}
\newtheorem*{proposition*}{Proposition}

\newtheorem*{example*}{Example}

\newtheorem*{conjecture*}{Conjecture}

\theoremstyle{definition}

\newtheorem*{definition*}{Definition}

\newtheorem*{notation*}{Notation}

\newtheorem{remark}[theorem]{Remark}
\newtheorem*{remark*}{Remark}

\theoremstyle{plain}

\newcommand{\theoremname}[1]{}

\newcommand{\gauge}[2]{{{#1}{.}{#2}}}
\newcommand{\imi}{ \bbi }

\newcommand{\Loop}{\Lambda}
\newcommand{\DOT}[2]{\langle#1,\,#2\rangle}

\DeclareMathOperator*{\spin}{spin}

\newcommand{\NEG}{\text{--}}

\newcommand{\Iso}{\mathrm{Iso}}

\newcommand{\arraystack}[2]{\begin{array}{c}#1\\#2\end{array}}

\newenvironment{smatrix}{\bigl[\begin{smallmatrix}}{\end{smallmatrix}\bigr]}

\newcommand{\tstack}[2]{\begin{tabular}{@{}c@{}}#1\\\footnotesize #2\end{tabular}}

\newcommand{\raiseimage}[1]{%
  \raisebox{-3ex}{\includegraphics{#1}}}

\newtoggle{usecolorimages}
\newcommand{\colorimages}{%
  \toggletrue{usecolorimages}}

\newcommand{\imagepath}{%
\iftoggle{usecolorimages}{image}{image-gray}}

\newcommand{\ximage}[4][0.3275]{%
  \tstack{\includegraphics[width=#1\textwidth]{\imagepath/#2}}%
  {#4%
  {\raisebox{-0.35cm}{\includegraphics[scale=0.9]{#3}}%
  }}}

\newgeometry{left=0.5in,right=2.6in,top=0.8in,bottom=0.8in}

\title[Cmc surfaces based on fundamental quadrilaterals]%
      {Constant mean curvature surfaces based on fundamental quadrilaterals}
\author{Alexander I. Bobenko}
\address{Institut f\"ur Mathematik, TU Berlin,
Str. des 17. Juni 136, 10623 Berlin, Germany}
\email{bobenko@math.tu-berlin.de}
\author{Sebastian Heller}
\address{Institut f\"ur Differentialgeometrie,
Universit\"at Hannover, 30167 Hannover, Germany}
\email{seb.heller@gmail.com}
\author{Nick Schmitt}
\address{Institut f\"ur Mathematik, TU Berlin,
Str. des 17. Juni 136, 10623 Berlin, Germany}
\email{schmitt@math.tu-berlin.de}
\date{\today}

\colorimages
\begin{document}

\begin{abstract}
We describe the construction of CMC surfaces with symmetries in
$\bbS^3$ and $\bbR^3$ using a CMC quadrilateral in a fundamental
tetrahedron of a tessellation of the space.  The fundamental piece is
constructed by the generalized Weierstrass representation using a
geometric flow on the space of potentials.
\end{abstract}
\maketitle
%\includefigure{r3-lattice3-a}

%\vspace{0.5cm}

\typeout{== section0 ============================================}\section*{Introduction}

Surfaces with constant mean curvature (CMC) in euclidean 3-space and
in the round 3-sphere can be investigated by methods of integrable
systems.  Their Gauss equation is the elliptic sinh-Gordon equation
\begin{equation}
\Delta u+ \sinh u=0\spacecomma
\end{equation}
which is one of the basic examples of integrable equations. Similar to
minimal surfaces in euclidean 3-space, CMC surfaces possess
1-parameter (denoted usually by $\lambda$) families of isometric
associated surfaces obtained by rotating their Hopf differential.
This allows CMC surfaces to be described in terms of loop
groups~\cite{Bobenko_1991:cmcsurfaces}, so that analytic methods of
the theory of integrable systems can be applied. One of the powerful
methods of the construction of CMC surfaces is the generalized
Weierstrass representation (DPW) by
Dorfmeister-Pedit-Wu~\cite{Dorfmeister_Pedit_Wu_1998}. It starts with
an analytic differential equation for the holomorphic frame $\Phi_z
=\Phi \xi$ with a meromorphic DPW potential $\xi(z,\lambda)$ and the
subsequent loop group factorization of $\Phi$, leading to immersion
formulas for the CMC surfaces. Control of the monodromy of the
holomorphic frame is of crucial importance for the construction of CMC
surfaces with non-trivial topology and symmetries.

\vspace{0.0975cm}%\begin{textblock*}{1cm}(1.5cm,1.5cm)
\begin{textblock*}{16cm}(14.5cm,7.5cm)
\begin{minipage}{5cm}
\begin{figure}
  \ximage[1]
    {r3-lattice3-cube-bulge.pdf}
    {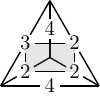}
    {\scriptsize (\textsc{a}) cube (bulge) \vspace{0.35cm}}
  \ximage[1]
    {r3-lattice3-cube-neck.pdf}
    {diagram/tet-i-r3-lattice3-cube.pdf}
    {\scriptsize (\textsc{b}) cube (neck) \vspace{0.35cm}}
  \ximage[1]
    {r3-lattice3-oct.pdf}
    {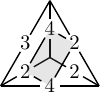}
    {\scriptsize (\textsc{c}) octahedron \vspace{0.35cm}}
  \captionsetup{margin=0cm}
  \caption{\label{fig:r3-lattice3-a}Triply periodic CMC surfaces in $\bbR^3$}
\end{figure}
\end{minipage}
\end{textblock*}

A particularly important class of potentials is given by Fuchsian
systems $\xi(z)$, those with only simple poles.  In the simplest case
of three singularities it reduces to the hypergeometric equation (see,
for example,~\cite{FokasItsKapaevNovokshenov}), whose monodromy group
can be described explicitly from the local residues. 
This leads to CMC surfaces based on
fundamental triangles~\cite{Schmitt_Kilian_Kobayashi_Rossman_2007}.
From the geometric point of view,  CMC surfaces
constructed from fundamental quadrilaterals are more natural, since they come from the
 curvature line parametrization.
 But for Fuchsian systems with
more then three singularities the monodromy cannot be computed
explicitly in terms of the coefficients of the system,
introducing accessory parameters.
 Then the simplest
holomorphic frame equation is a Fuchsian system with four
singularities on the Riemann sphere
\begin{equation}
\label{eq:Fuchsian4}
\Phi_z=\Phi \sum_{k=0}^3 \frac{A_k}{z-z_k}\spaceperiod
\end{equation}
In \cref{sec:FuchsianDPW4} of this paper we show how all periodic and
compact surfaces based on fundamental quadrilaterals can be
constructed from the system \eqref{eq:Fuchsian4}.  Our constructions
make explicit use of this Fuchsian DPW form. The relation of the
monodromy and the coefficients of the Fuchsian system is the famous
Hilbert's 21st problem, which was intensely
studied~\cite{AnosovBolibruch}. 
 There exist
many important partial results in the simplest non-trivial case of
four singularities. 
This case was investigated mostly within the
theory of isomonodromic deformations~\cite{FokasItsKapaevNovokshenov}
and the Painlev\'{e} VI equation, where the problem is to describe the
coefficients $A_k$ as functions of the poles $z_j$ when the monodromy
group is preserved. The holomorphic frame $\Phi(z,\lambda)$ of a CMC
surface lies in a loop group, and the main analytic problem is to
construct solutions whose monodromy group is unitary on the unit
circle $|\lambda |=1$, giving global solutions of the Gauss equation
on the four-punctured sphere.
%\includefigure{r3-lattice3-a}

In general, it is a hard problem to control the intrinsic and
extrinsic closing conditions to obtain closed surfaces or surfaces
with prescribed global properties.  In recent years, important
progress has been made using a flow of DPW
potentials~\cite{Traizet_2020_Crelle, Traizet_20_MA,
  Heller_Heller_Traizet_19} or similar methods on spectral
data~\cite{Heller_Heller_Schmitt_2018}.  By the very nature of these
techniques, only surfaces which are small perturbations of spheres or
tori have been reached~\cite{Heller_Heller_Traizet_19}.

\vspace{0.0975cm}\begin{textblock*}{16cm}(14.5cm,5cm)
\begin{minipage}{5cm}
\begin{figure}
  \ximage[1]
    {r3-lattice3-altoct.pdf}
    {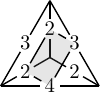}
    {\scriptsize (\textsc{a}) alternate octahedron \vspace{0.35cm}}
  \ximage[1]
    {r3-lattice3-altoct-b.pdf}
    {diagram/tet-i-r3-lattice3-altoct.pdf}
    {\scriptsize (\textsc{b}) alternate octahedron \vspace{0.35cm}}
  \ximage[1]
    {r3-lattice3-altcube-c.pdf}
    {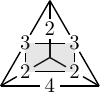}
    {\scriptsize (\textsc{c}) alternate cube \vspace{0.35cm}}
  \captionsetup{margin=0cm}
  \caption{\label{fig:r3-lattice3-b}Triply periodic CMC surfaces in $\bbR^3$}
\end{figure}
\end{minipage}
\end{textblock*}

In~\cite{Lawson_1970} Lawson constructed the first compact minimal
surfaces in the round 3-sphere of genus $g\geq2.$ A fundamental piece
of a Lawson surface is obtained by the Plateau solution of a specific
geodesic polygon. The compact surface is then built from the
fundamental 4-gon by the finite group generated by rotations around
the geodesic edges of the polygon. Later,
Karcher-Pinkall-Sterling~\cite{Karcher_Pinkall_Sterling_1988}
constructed new minimal surfaces in the 3-sphere by starting with a
tessellation of the 3-sphere into tetrahedra.  The minimal surfaces
are obtained from fundamental minimal 4-gons within such a tetrahedron
which reflect across the geodesic boundaries. Constant mean curvature
(CMC) surfaces in $\bbR^3$ have been constructed by adapting these
methods~\cite{Grosse-Brauckmann_1993}; see also~\cite{KPOb,KGBKP} for
related computer experiments.

 This paper constructs such fundamental patches of surfaces
 (\cref{sec:experimental-part}) based on the deformation of DPW
 potentials.  In this paper the following new surfaces are numerically
 constructed: triply periodic surfaces
 (\cref{fig:r3-lattice3-a}\textsc{b-c}, \cref{fig:r3-lattice3-b}) and
 doubly periodic surfaces (\cref{fig:r3-lattice2-a}\textsc{b-c},
 \cref{fig:r3-lattice2-b}\textsc{b-c},
 \cref{fig:r3-lattice2-c}\textsc{b-c}), new doubly periodic surfaces
 with Delaunay ends (\cref{fig:r3-lattice2-ends}), new surfaces with
 Delaunay ends of positive genus
 (\cref{fig:r3-platonic-a,fig:r3-platonic-b,fig:s3-torus}) as well as
 new KPS-type surfaces (\cref{fig:s3-platonic-b}\textsc{a,d}).  We
 also reconstruct by these methods previously constructed surfaces
 based on doubly-periodic hexagonal, square and triangular tilings of
 the plane, triply periodic cubic examples, cylinders with
 ends~\cite{Grosse-Brauckmann_1993,KGB2}, and the Lawson and KPS
 surfaces~\cite{Lawson_1970,Karcher_Pinkall_Sterling_1988}.

The 3D-data of the surfaces constructed in this paper are available in
the DGD Gallery~\cite{DGD-Gallery}.

\section*{Acknowledgements}
The first author is partially supported by the DFG Collaborative
Research Center TRR 109 \emph{Discretization in Geometry and
Dynamics}.  The second author is supported by the DFG grant HE
6829/3-1 of the DFG priority program SPP 2026 \emph{Geometry at
Infinity}.  The third author is supported by the DFG Collaborative
Research Center TRR 109 \emph{Discretization in Geometry and
Dynamics}.

\typeout{== section1 ============================================}\section{Geometric construction}
\label{sec:geometric_construction}

\subsection{The construction}
\label{sec:overview}

\begin{wrapfigure}{r}{0.1\linewidth}
\includegraphics[scale=1.25]{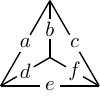}
\end{wrapfigure}
This report reports on the experimental construction of CMC surfaces
in $\bbR^3$ with genus with and without Delaunay ends via the
generalized Weierstrass representation
(DPW)~\cite{Dorfmeister_Pedit_Wu_1998}.  The construction starts with
a tetrahedron in $\bbR^3$ as shown which tessellates $\bbR^3$ by the
group generated by the reflections in the four planes containing its
faces.  Each of the six edges of the tetrahedron is marked with an
integer $n\in\bbN_{\ge 1}\cup\{\infty\}$ specifying that the internal
dihedral angle between the two planes meeting at that edge is $\pi/n$.
The tetrahedron can be degenerate in the following ways:

\begin{itemize}
\item vertex at $\infty$
\item parallel planes, with opposite outward normals: the edge of the
  tetrahedron between the two planes is marked with $\infty$
\item coincident planes, with the same outward normal: the edge of the
  tetrahedron between the two planes is marked with $1$.
\end{itemize}

\begin{wrapfigure}{r}{0.1\linewidth}
\includegraphics[scale=1.25]{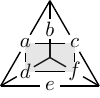}
\end{wrapfigure}
In this tetrahedron construct a CMC quadrilateral as shown, such that
\begin{itemize}
\item
  each of the four edges of the quadrilateral lies in a plane of the
  tetrahedron, and the surface reflects smoothly across this plane
\item
  at each of the four vertices of the quadrilateral, application of
  the tessellation group results in a surface with either an immersed
  point of the surface or a once-wrapped Delaunay end at the vertex.
  \end{itemize}
Then the surface constructed by application of the tessellation group
is a CMC immersion with optional Delaunay ends and the symmetries of
that group.  Its genus is finite if the tessellation group is finite,
and infinite if the group is infinite.  These surfaces are described
in detail at the end of this section.

The polygon is constructed via a Fuchsian DPW potential on $\CPone$
with four simple poles on $\bbS^1$ (\cref{sec:potential}) and a
reflection symmetry across $\bbS^1$.  The unit disk is the domain of a
CMC quadrilateral which reflects in planes containing its boundaries.
The simple poles with constant or Delaunay residue eigenvalues insure
that each vertex of the quadrilateral after reflection is either
immersed or a Delaunay end.  The four dihedral angles of the
tetrahedron at the corners of the quadrilateral are controlled by the
four local monodromies of the potential, and the two remaining
dihedral angles by two global monodromies.

The potential has two accessory parameters which are computed by the
unitary flow (\cref{sec:unitary-flow}).  Starting with an initial
surface (\cref{sec:initial}) which satisfies the intrinsic closing
condition (unitary monodromy on $\bbS^1_\lambda$), the unitary flow,
which preserves this condition, is run through the space of potentials
until the dihedral angles of the planes reach the values prescribed by
the tessellation.  The dihedral angles are controlled by certain
monodromy traces at the evaluation point.  The unitary flow is not
known in general to exist, but short time existence can be shown in
some cases~\cite{Heller_Heller_Schmitt_2018}. Hence we construct the
surfaces numerically, giving evidence that the unitary flow has long
time existence.

The Lawson surfaces~\cite{Lawson_1970} (\cref{fig:s3-lawson}) and the
surfaces of
Karcher-Pinkall-Sterling~\cite{Karcher_Pinkall_Sterling_1988}
(\cref{fig:s3-platonic-a,fig:s3-platonic-b} and
\cref{fig:s3-cell-a,fig:s3-cell-b}) have been constructed by solutions
of Plateau problems.  The cubic lattice and some of the 2-dimensional
lattices were shown to exist by similar
methods~\cite{Grosse-Brauckmann_1993}.

\subsection{Tetrahedral tessellations}

The following theorem classifies the tetrahedral tessellation of
$\bbS^3$ (which are compact) and of $\bbR^3$ (which are compact,
paracompact or degenerate). The tetrahedral tessellation of $\bbH^3$,
which can be determined by the same methods, are omitted for
simplicity.

\vspace{0.0975cm}\begin{textblock*}{16cm}(14.5cm,7.5cm)
\begin{minipage}{5cm}
\begin{figure}
  \ximage[1]
    {r3-lattice2-hexagon-bulge.pdf}
    {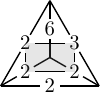}
    {\scriptsize (\textsc{a}) hexagon (bulge) \vspace{0.35cm}}
  \ximage[1]
    {r3-lattice2-hexagon-neck.pdf}
    {diagram/tet-i-r3-lattice2-hexagon.pdf}
    {\scriptsize (\textsc{b}) hexagon (neck) \vspace{0.35cm}}
  \ximage[1]
    {r3-lattice2-althexagon.pdf}
    {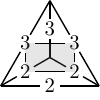}
    {\scriptsize (\textsc{c}) alternate hexagon \vspace{0.35cm}}
  \captionsetup{margin=0cm}
  \caption{\label{fig:r3-lattice2-a}Doubly periodic CMC surfaces in $\bbR^3$}
\end{figure}
\end{minipage}
\end{textblock*}

\begin{theorem}
  \label{thm:tessellation}
  \mbox{}

\begin{enumerate}
\item
  \label{item:tessellation-1}
  The tetrahedral tessellations of $\bbS^3$ are as follows:
  \begin{equation*}
    \begin{tabular}{>{\footnotesize}l|>{\footnotesize}l}
          \raiseimage{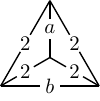}
      \hspace{0.5em} $a,\,b\in\bbN_{\ge 2}$ &
      \parbox{0.4\textwidth}{Minimal Lawson surfaces $\xi_{ab}$}
      \\ \raiseimage{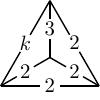}
      \hspace{0.5em} $k\in\{3,\,4,\,5\}$ &
      \parbox{0.4\textwidth}{ surfaces with respective tetrahedral,
        octahedral, and icosahedral symmetries}
      \\ \raiseimage{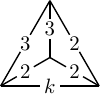}
      \hspace{0.5em} $k\in\{3,\,4,\,5\}$ &
      \parbox{0.4\textwidth}{ surfaces with respective $5$-cell,
        $16$/$8$-cell and $600$/$120$-cell symmetries}
      \\ \raiseimage{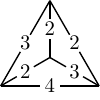}
      \hspace{0.5em} &
      \parbox{0.4\textwidth}{surface with $24$-cell symmetry}
      \\ \raiseimage{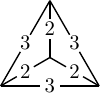}
      \hspace{0.5em} &
      \parbox{0.25\textwidth}{subgroup of $16$-cell}
    \end{tabular}
  \end{equation*}
\item
  \label{item:tessellation-2}
  The tetrahedral tessellations of $\bbR^3$ are as follows:
  \begin{equation*}
    \begin{tabular}{>{\footnotesize}l|>{\footnotesize}l}
      \raiseimage{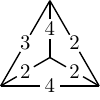}
      \hspace{0.5em} \raiseimage{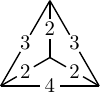}
      \hspace{0.5em} \raiseimage{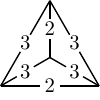} &
      \parbox{0.4\textwidth}{triply periodic surfaces}
      \\ \raiseimage{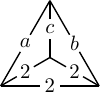}
      \hspace{0.5em}
      \parbox{0.25\textwidth}{$(a,\,b,\,c)$ one of $(3,\,3,\,3)$,
        $(2,\,3,\,6)$ or $(2,\,4,\,4)$} &
      \parbox{0.4\textwidth}{doubly periodic surfaces}
      \\ \raiseimage{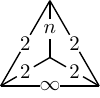}
      \hspace{0.5em}
      \parbox{0.25\textwidth}{$n\in\bbN_{\ge 2}\cup\{\infty\}$} &
      \parbox{0.4\textwidth}{singly periodic surfaces}
      \\ \raiseimage{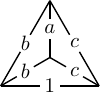}
      \hspace{0.5em}
      \parbox{0.25\textwidth}{$(a,\,b,\,c)$ a permutation of
        $(2,\,3,\,3)$, $(2,\,3,\,4)$ or $(2,\,3,\,5)$} &
      \parbox{0.4\textwidth}{ surfaces with Platonic symmetries and
        Delaunay ends}
      \\ \raiseimage{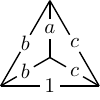}
      \hspace{0.5em}
      \parbox{0.25\textwidth}{$(a,\,b,\,c)$ a permutation of
        $(3,\,3,\,3)$, $(2,\,4,\,4)$ or $(2,\,3,\,6)$}
      \\ \raiseimage{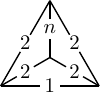}
      \hspace{0.5em}
      \raiseimage{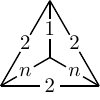}
      \hspace{0.5em} $n\in\bbN_{\ge 2}\cup\{\infty\}$ &
      \parbox{0.4\textwidth}{tori with Delaunay ends}
      \\ \raiseimage{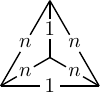}
      \hspace{0.5em}
      \raiseimage{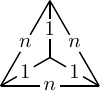}
      \hspace{0.5em} $n\in\bbN_{\ge 2}\cup\{\infty\}$
      \\ \raiseimage{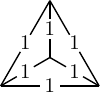}
    \end{tabular}
    \end{equation*}
\end{enumerate}
\end{theorem}

\begin{proof}
  Necessary conditions that a compact tetrahedron tessellates one of
  the spaceforms $\bbS^3$, $\bbR^3$ or $\bbH^3$ are the following.
  \begin{itemize}
  \item
    Each edge of the tetrahedron is marked with an integer
    $n\in\bbN_{\ge 2}$ denoting that the internal dihedral angle
    between the two faces meeting at that edge is $\pi/n$.
  \item At each vertex of the tetrahedron, the three integers marking
    the three edges meeting at the vertex are $(2,\,2,\,n)$,
    $n\in\bbN_{\ge 2}$ or $(2,\,3,\,k)$, $k\in\{3,\,4,\,5\}$.
  \item
    The Gram matrix $T\in\matN{n}(\bbR)$ defined by $T_{ij} =
    -\cos\pi/n_{ij}$ has signature $(\delta,\,1,\,1,\,1)$, where
    $\delta=1$, $0$ or $-1$ for $\bbS^3$, $\bbR^3$ and $\bbH^3$
    respectively.
  \end{itemize}

\vspace{0.0975cm}\begin{textblock*}{16cm}(14.5cm,7.5cm)
\begin{minipage}{5cm}
\begin{figure}
  \ximage[1]
    {r3-lattice2-square-bulge.pdf}
    {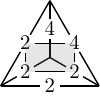}
    {\scriptsize (\textsc{a}) square (bulge) \vspace{0.35cm}}
  \ximage[1]
    {r3-lattice2-square-neck.pdf}
    {diagram/tet-i-r3-lattice2-square.pdf}
    {\scriptsize (\textsc{b}) square (neck) \vspace{0.35cm}}
  \ximage[1]
    {r3-lattice2-altsquare.pdf}
    {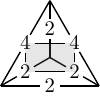}
    {\scriptsize (\textsc{c}) alternate square \vspace{0.35cm}}
  \captionsetup{margin=0cm}
  \caption{\label{fig:r3-lattice2-b}Doubly periodic CMC surfaces in $\bbR^3$}
\end{figure}
\end{minipage}
\end{textblock*}

  The compact tetrahedra in $\bbS^3$, $\bbR^3$ and $\bbH^3$ are the
  following:
  \begin{equation*}
    \arraystack{\includegraphics{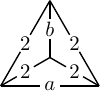}}{
      \text{\tiny $a,b\in\bbN_{\ge 2}$}\\A_{ab}} \quad
    \arraystack{\includegraphics{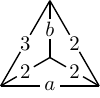}}{
      \text{\tiny $a,b\in\{2,3,4,5\}$}\\B_{ab}} \quad
    \arraystack{\includegraphics{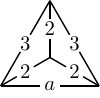}}{
      \text{\tiny $a\in\{2,3,4,5\}$}\\C_{a}} \quad
    \arraystack{\includegraphics{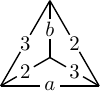}}{
      \text{\tiny $a,b\in\{2,3,4,5\}$}\\D_{ab}}
\end{equation*}
with identifications $A_{ba} = A_{ab}$, $B_{ba} = B_{ab}$, $D_{ba} =
D_{ab}$ and
\begin{equation}
  \label{eq:tet-identities}
A_{23} = B_{22}\spacecomma\quad A_{33} = D_{22}\spacecomma\quad B_{23}
= C_{2}\spacecomma\quad B_{33} = D_{23}\spaceperiod
\end{equation}
To see this, first consider those tetrahedra with at least one edge
marked with $4$ or $5$.  Then at each of the vertices at the endpoints
of that edge, the other two edges meeting the vertex must be marked
with $2$ and $3$.  Hence all such tetrahedra with at least one $4$ or
$5$ is one of the four types $A_{ab}$, $B_{ab}$, $C_{a}$ or $D_{ab}$.
The remaining tetrahedra have only $2$ or $3$ at each face.  There are
seven of these, namely $A_{22}$, $A_{23}=B_{22}$, $B_{23}=C_{2}$,
$C_{3}$, $A_{33}=D_{22}$, $B_{33}=D_{23}$ and $D_{33}$.

Since its Gram matrix has positive determinant, the tetrahedron
$A_{ab}$ is in $\bbS^3$.  The spaceforms for the other tetrahedra
$B_{ab}$, $C_{a}$ and $D_{ab}$ are determined by the signs of the
determinate of the Gram matrix as follows:
\begin{equation}
  \small
  \begin{array}[t]{c|cccc}
    B & 2 & 3 & 4 & 5\\ \hline 2 & [+] & + & + & +\\ 3 & & + & + &
    +\\ 4 & & & 0 & -\\ 5 & & & & -
  \end{array}
  \qquad
  \begin{array}[t]{c|cccc}
    C & 2 & 3 & 4 & 5\\ \hline & [+] & + & 0 & -
  \end{array}
  \qquad
  \begin{array}[t]{c|cccc}
    D & 2 & 3 & 4 & 5\\ \hline 2 & [+] & [+] & + & -\\ 3 & & 0 & - &
    -\\ 4 & & & - & -\\ 5 & & & & -
  \end{array}
\end{equation}
The $[+]$ in the above tables denotes entries which are redundant due
to the identifications~\eqref{eq:tet-identities}

Hence the tetrahedra which tessellate $\bbS^3$ are $A_{ab}$ and the
seven tetrahedra $B_{23}$, $B_{24}$, $B_{25}$, $B_{33}$, $B_{34}$,
$B_{35}$ and $D_{24}$. These tessellate $\bbS^3$ because they
tessellate either a sphere or a $n$-cell, which in turn tessellates
$\bbS^3$.

The compact tetrahedra which tessellate $\bbR^3$ are the three
tetrahedra $B_{44}$, $C_{4}$ and $D_{33}$.  The first two tessellate a
cube and the third tessellates a rhombic dodecahedron, each of which
in turn tessellates $\bbR^3$.

The paracompact tetrahedral tessellations of $\bbR^3$ are classified
similarly except that the integer triple at each vertex is as in the
compact case or one of $(3,\,3,\,3)$, $(2,\,4,\,4)$ or $(2,\,3,\,6)$.

The degenerate tetrahedral tessellations of $\bbR^3$ are classified
similarly except that the integer triple at each vertex is as in the
paracompact case or one of $(2,\,2,\,\infty)$ or $(1,\,n,\,n)$,
$n\in\bbN_{\ge 1}\cup\{\infty\}$.
\end{proof}

\subsection{The surfaces}
\label{sec:surface}

This section describes the experimentally constructed minimal surfaces
in $\bbS^3$ and CMC surfaces in $\bbR^3$. They are of the types:
\newline
In $\bbR^3$:

\begin{itemize}
\item
  triply periodic CMC surfaces $\bbR^3$ without ends
  (\cref{fig:r3-lattice3-a,fig:r3-lattice3-b})
\item
  doubly periodic CMC surfaces $\bbR^3$ without ends
  (e.g. \cref{fig:r3-lattice2-a})
\item
  doubly periodic CMC surfaces $\bbR^3$ with ends
  (\cref{fig:r3-lattice2-ends})
\item
  single periodic CMC surfaces in $\bbR^3$ with Delaunay ends
  (cylinders, \cref{fig:r3-cylinder})
\item
  CMC surfaces in $\bbR^3$ with dihedral symmetry and Delaunay ends
  (tori, \cref{fig:r3-torus-a,fig:r3-torus-b})
\item
  CMC surfaces in $\bbR^3$ with Platonic symmetries and Delaunay ends
  (\cref{fig:r3-platonic-a,fig:r3-platonic-b})
\item
  CMC spheres in $\bbR^3$ with four Delaunay ends (fournoids,
  \cref{fig:r3-fournoid}).
\end{itemize}

In $\bbS^3$:
\begin{itemize}
\item
  Lawson surfaces $\xi_{ab}$ in $\bbS^3$ (\cref{fig:s3-lawson})
\item
  minimal surfaces in $\bbS^3$ with Platonic symmetries
  (\cref{fig:s3-platonic-a,fig:s3-platonic-b})
\item
  minimal surfaces in $\bbS^3$ with $n$-cell symmetries
  (\cref{fig:s3-cell-a,fig:s3-cell-b})

\item
  minimal tori in $\bbS^3$ with Delaunay ends (\cref{fig:s3-torus}).
\end{itemize}
In all shown figures the lines on the surfaces are curvature lines.
The surfaces in $\bbS^3$ are stereographically projected to $\bbR^3$.

\vspace{0.0975cm}\begin{textblock*}{16cm}(14.5cm,7.5cm)
\begin{minipage}{5cm}
\begin{figure}
  \ximage[1]
    {r3-lattice2-triangle-bulge.pdf}
    {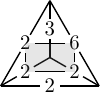}
    {\scriptsize (\textsc{a}) triangle (bulge) \vspace{0.35cm}}
  \ximage[1]
    {r3-lattice2-triangle-neck.pdf}
    {diagram/tet-i-r3-lattice2-triangle.pdf}
    {\scriptsize (\textsc{b}) triangle (neck) \vspace{0.35cm}}
  \ximage[1]
    {r3-lattice2-rhombus.pdf}
    {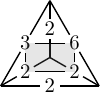}
    {\scriptsize (\textsc{c}) rhombus \vspace{0.35cm}}
  \captionsetup{margin=0cm}
  \caption{\label{fig:r3-lattice2-c}Doubly periodic CMC surfaces in $\bbR^3$}
\end{figure}
\end{minipage}
\end{textblock*}

\textbf{Triply periodic surfaces in $\bbR^3$.}

\begin{figure}[H]
\includegraphics[height=0.1\textheight]{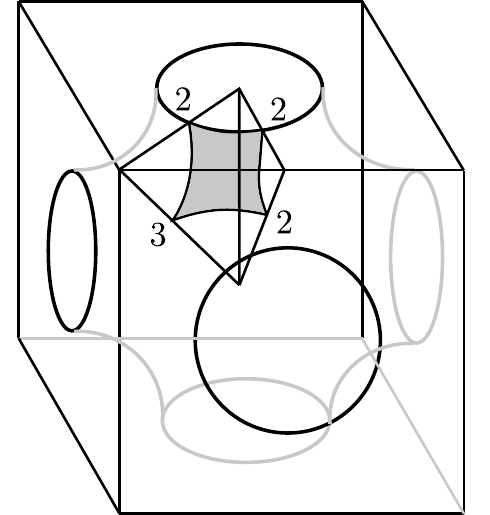}
\caption{The simplest triply periodic surface and its tetrahedron.}
\label{fig:diagram-cube}
\end{figure}

The simplest triply periodic surfaces can be thought of as tubes along
the edges of the standard cubic lattice in $\bbR^3$
(\cref{fig:diagram-cube}).  The genus of those surfaces modulo
translation is $3$.

\begin{figure}[H]
  \begin{tabular}{l|l|l}
    \tstack{\includegraphics[scale=1]{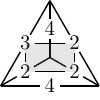}}{cube
      $a_1$}
    \hspace{0.5em}
    \tstack{\includegraphics[scale=1]{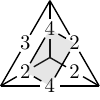}}{oct
      $a_2$}
    \hspace{0.5em}
    \tstack{\includegraphics[scale=1]{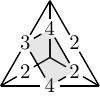}}{$a_3$}
    \quad & \quad
    \tstack{\includegraphics[scale=1]{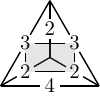}}{alt
      cube $b_1$}
  \hspace{0.5em}
  \tstack{\includegraphics[scale=1]{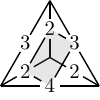}}{alt
    oct $b_2$} \quad & \quad
  \tstack{\includegraphics[scale=1]{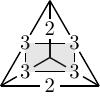}}{$c_1$}
  \hspace{0.5em}
  \tstack{\includegraphics[scale=1]{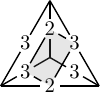}}{$c_2$}
  \end{tabular}
  \caption{The tetrahedra for the seven possible triply periodic
    examples.}
  \label{fig:tetR37}
\end{figure}

The triply periodic surfaces are constructed from the three compact
tetrahedra which tessellate $\bbR^3$.  Since the CMC quadrilateral can
be situated in each tetrahedron in three ways, this gives nine
different configurations, of which two are redundant due to symmetry
(\cref{fig:tetR37}).
Of these, $a_3$ and $b_2$ are not possible under our symmetry
constraints (compare with \eqref{symmetric-potential}), $c_1$ seems to
devolve to $a_1$, and the flow for $c_2$ degenerates
(\cref{fig:r3-lattice3-a,fig:r3-lattice3-b}).

\textbf{Doubly periodic surfaces in $\bbR^3$.}

The doubly periodic surfaces can be thought of as tubes along the
edges of a triangle tessellations of $\bbR^2$.  The six 2-dimensional
lattices are constructed with the tetrahedron below where
$(a,\,b,\,c)$ are the indices of a triangle tessellation of $\bbR^2$,
that is, a permutation of $(3,\,3,\,3)$, $(2,\,4,\,4)$ or
$(2,\,3,\,6)$
(\cref{fig:r3-lattice2-a,fig:r3-lattice2-b,fig:r3-lattice2-c}).

\begin{figure}[H]
  \raisebox{-4ex}{\includegraphics[scale=1.5]{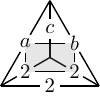}}
  \qquad\qquad
  \begin{tabular}{>{\footnotesize}l|>{\footnotesize}r|>{\footnotesize}r}
    lattice & (a,\,b,\,c) & genus \\ \hline hexagon & (2,\,3,\,6) & 2
    \\ square & (2,\,3,\,6) & 2 \\ triangle & (2,\,3,\,6) & 3
    \\ \hline alt hexagon & (3,\,3,\,3) & 2 \\ alt square &
    (4,\,4,\,2) & 2 \\ rhombus & (3,\,6,\,2) & 2
\end{tabular}
  \qquad \qquad
  \raisebox{-1.0cm}{\includegraphics[height=0.08\textheight]{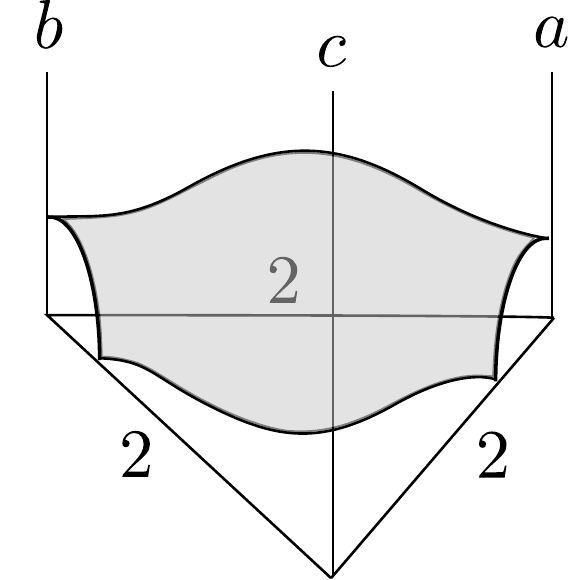}
  }
\caption{ Left: tetrahedron for doubly periodic surface,, where
  $(a,\,b,\,c)$ is a permutation of $(3,\,3,\,3)$, $(2,\,3,\,6)$ or
  $(2,\,4,\,4)$.  Middle: table of 2-dimensional lattices.  The genus
  listed in the table is that of the surface modulo translations.
  Right: fundamental piece for a 2-dimensional lattice.}
\end{figure}

\begin{figure}[H]
  \includegraphics[width=0.3\textwidth]{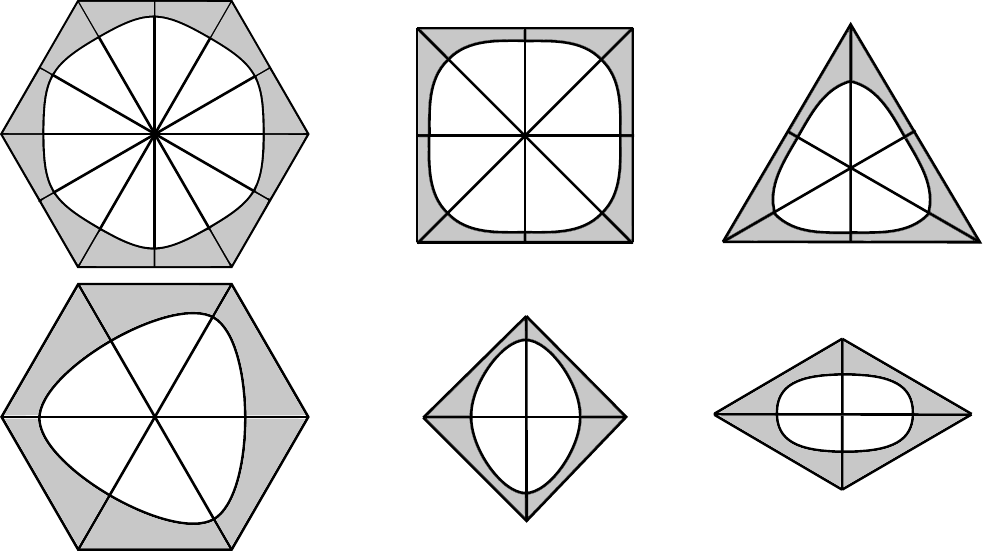}
  \caption{The six $2$-dimensional lattices.}
\end{figure}

\textbf{Doubly periodic surfaces in $\bbR^3$ with Delaunay ends.}  The
doubly periodic surfaces with Delaunay ends are obtained from triangle
tessellations of $\bbR^2$.  Additional freedom is given by the choice
of vertices corresponding to Delaunay ends
(\cref{fig:r3-lattice2-ends}).

\textbf{Cylinders in $\bbR^3$ with ends.}  Cylinders with Delaunay
ends can be constructed from a degenerate tetrahedron with two
parallel planes.  Of course, the same construction without Delaunay
ends give the classical rotational symmetric periodic surfaces, i.e.,
Delaunay cylinders (\cref{fig:r3-cylinder}).

\textbf{Tori in $\bbR^3$ with Delaunay ends.}  The torus with $n$ ends
is constructed via the diagram below
(\cref{fig:r3-torus-a,fig:r3-torus-b}).  For large $n$, existence of
those tori can be shown by growing Delaunay ends in equidistance on
one side of a cylinder.

\begin{figure}[H]
  \raisebox{4ex}{\includegraphics[scale=1.5]{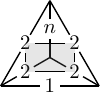}}
  \qquad\qquad\qquad
  \includegraphics[height=0.12\textheight]{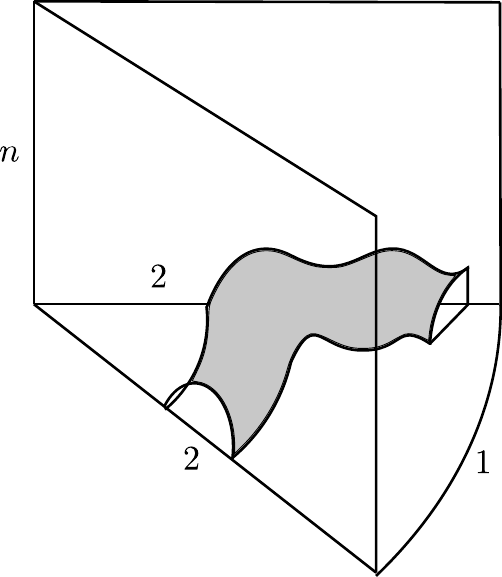}
  \caption{Left: tetrahedron for a torus in $\bbR^3$ with $n$ Delaunay
    ends and order $n$ cyclic symmetry and $n$ ends. Right:
    fundamental piece of this torus.}
\end{figure}

\vspace{0.0975cm}\begin{textblock*}{16cm}(14.5cm,13.8cm)
\begin{minipage}{5cm}
\begin{figure}
  \ximage[1]
    {r3-lattice2-ends-hexagon-bulge.pdf}
    {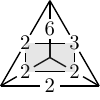}
    {\scriptsize (\textsc{a}) hexagon lattice with ends (bulge) \vspace{0.35cm}}
  \ximage[1]
    {r3-lattice2-ends-hexagon-neck.pdf}
    {diagram/tet-i-r3-lattice2-ends-hexagon.pdf}
    {\scriptsize (\textsc{b}) hexagon lattice with ends (neck) \vspace{0.35cm}}
  \captionsetup{margin=0cm}
  \caption{\label{fig:r3-lattice2-ends}Doubly periodic CMC surfaces in $\bbR^3$ with Delaunay ends}
\end{figure}
\end{minipage}
\end{textblock*}

\textbf{Surfaces in $\bbR^3$ with Platonic symmetry and Delaunay
  ends.}

Given a triangle tessellations of $\bbS^2$, the surface is the orbit
of a tube along one edge of the triangle with a Delaunay end at a
vertex of the triangle.  Equivalently, the surface is built from tubes
along the edges of one of the five Platonic solids, with ends
emanating from the vertices.  The five tetrahedra are as in the
diagram below, with $(a,\,b,\,c)$ a permutation of $(2,\,3,\,k)$,
$k\in\{3,\,4,\,5\}$ (\cref{fig:r3-platonic-a,fig:r3-platonic-b}).

\begin{figure}[H]
  \raisebox{0.5cm}{
    \raisebox{-8ex}{\includegraphics[scale=1.5]{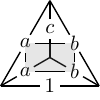}}}
  \qquad
  \begin{tabular}{>{\footnotesize}l|>{\footnotesize}r|>{\footnotesize}r}
    surface & (a,\,b,\,c) & genus \\ \hline tetrahedron & (2,\,3,\,6)
    & 3 \\ octahedron & (2,\,3,\,6) & 7 \\ icosahedron & (2,\,3,\,6) &
    19 \\ \hline cube & (2,\,3,\,6) & 5 \\ dodecahedron & (2,\,3,\,6)
    & 11
  \end{tabular}
  \qquad \raisebox{-1.4cm}{%
    \includegraphics[height=0.12\textheight]{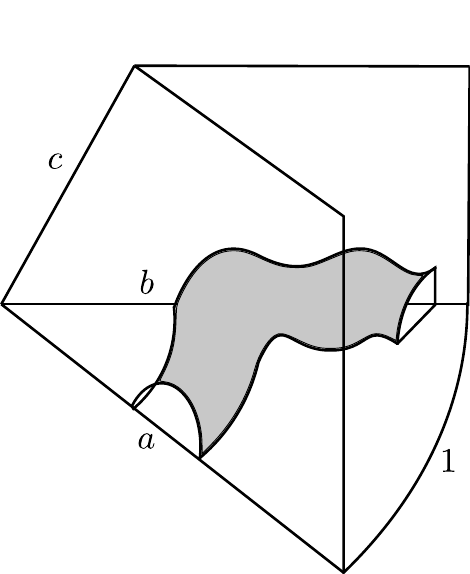} }
\caption{Left: Tetrahedron for Platonic surfaces in $\bbR^3$, where
  $(a,\,b,\,c)$ is a permutation of $(2,\,3,\,k)$.  Right: fundamental
  piece of the surfaces with Platonic symmetry and Delaunay ends.}
\end{figure}

\textbf{Lawson surfaces.}  Classically, the Lawson
surfaces~\cite{Lawson_1970} are constructed from Plateau solutions of
a geodesic polygon by reflection.  The tetrahedron $A_{ab}$ and its
inscribed fundamental piece admit a rotational order 2 symmetry around
a geodesic through the vertices labeled by $a$ and $b$.  The geodesic
arc is contained in the fundamental piece. This observation relates
the original construction with the construction carried out in the
present work, see also~\cite{Karcher_Pinkall_Sterling_1988}
(\cref{fig:s3-lawson}).

\begin{figure}[H]
  \includegraphics[scale=1.5]{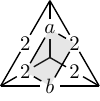}
  \caption{Tetrahedron for Lawson surface $\xi_{a-1,b-1}$.}
\end{figure}

\textbf{Surfaces in $\bbS^3$ with Platonic symmetries.}

Minimal surfaces in $\bbS^3$ with Platonic symmetries have been
constructed by
Karcher-Pinkall-Sterling~\cite{Karcher_Pinkall_Sterling_1988}.  These
surfaces can be thought of as tubes along one edge of a triangle which
tessellates $\bbS^2$. Note that~\cite{Karcher_Pinkall_Sterling_1988}
does not list all possible surfaces, e.g. the alternate octahedron of
genus 11 and the alternate icosahedron of genus 29 are missing
(\cref{fig:s3-platonic-a,fig:s3-platonic-b}).

\begin{figure}[H]
  \raisebox{-5ex}{\includegraphics[scale=1.5]{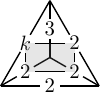}}
  \qquad\qquad
  \begin{tabular}{>{\footnotesize}l|>{\footnotesize}l}
    symmetry & genus\\ \hline tetrahedron & 3\\ cube & 5\\ octahedron
    & 7\\ alt octahedron & 11 \\ dodecahedron & 11\\ icosahedron &
    19\\ alt icosahedron & 29
  \end{tabular}
  \caption{Left: tetrahedron for surfaces in $\bbS^3$ with Platonic
    symmetries.  Right: table of these surfaces.}
\end{figure}

\begin{figure}[H]
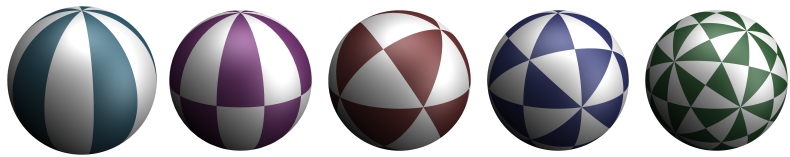

  \includegraphics[width=0.15\textwidth]{\imagepath/s2tile-cyc.pdf}
  \includegraphics[width=0.15\textwidth]{\imagepath/s2tile-dih.pdf}
  \includegraphics[width=0.15\textwidth]{\imagepath/s2tile-tet.pdf}
  \includegraphics[width=0.15\textwidth]{\imagepath/s2tile-oct.pdf}
  \includegraphics[width=0.15\textwidth]{\imagepath/s2tile-ico.pdf}
  \caption{Triangle tessellations of $\bbS^2$.  Left to right: cyclic
    of order $5$, dihedral of order $10$, tetrahedral, octahedral and
    icosahedra.}
\end{figure}

\vspace{0.0975cm}\begin{textblock*}{16cm}(14.5cm,7.5cm)
\begin{minipage}{5cm}
\begin{figure}
  \ximage[1]
    {r3-cylinder-2.pdf}
    {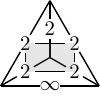}
    {\scriptsize (\textsc{a}) cylinder 2 \vspace{0.35cm}}
  \ximage[1]
    {r3-cylinder-3.pdf}
    {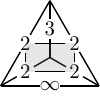}
    {\scriptsize (\textsc{b}) cylinder 3 \vspace{0.35cm}}
  \ximage[1]
    {r3-cylinder-4.pdf}
    {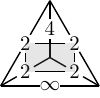}
    {\scriptsize (\textsc{c}) cylinder 4 \vspace{0.35cm}}
  \captionsetup{margin=0cm}
  \caption{\label{fig:r3-cylinder}CMC cylinders in $\bbR^3$ with Delaunay ends}
\end{figure}
\end{minipage}
\end{textblock*}

\textbf{Surfaces in $\bbS^3$ with $n$-cell symmetries.}

For each of the $n$-cell tessellations of $\bbS^3$ there is a surface
which can be thought of as tubes along the edges of the cells. These
minimal surfaces have also been constructed by
Karcher-Pinkall-Sterling~\cite{Karcher_Pinkall_Sterling_1988}
(\cref{fig:s3-cell-a,fig:s3-cell-b}).

\begin{figure}[H]
  \raisebox{-5ex}{\includegraphics[scale=1.5]{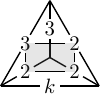}}
  \qquad\qquad
  \begin{tabular}{>{\footnotesize}l|>{\footnotesize}l}
    \text{symmetry} & \text{genus}\\ \hline $5$-cell & 6\\ $16$-cell &
    17\\ $24$-cell & 73\\ $600$-cell & 601
  \end{tabular}
  \caption{Left: tetrahedron for surfaces in $\bbS^3$ with $n$-cell
    symmetries.  Right: table of these surfaces.}
\end{figure}

\begin{figure}[H]
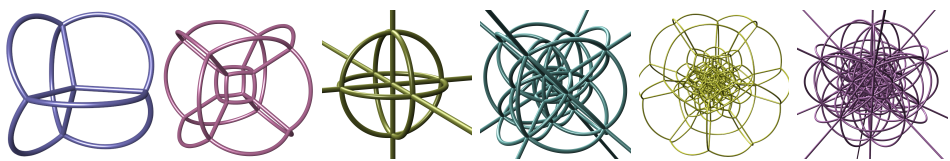

  \includegraphics[width=0.15\textwidth]{\imagepath/cell-5.pdf}
  \includegraphics[width=0.15\textwidth]{\imagepath/cell-8.pdf}
  \includegraphics[width=0.15\textwidth]{\imagepath/cell-16.pdf}
  \includegraphics[width=0.15\textwidth]{\imagepath/cell-24.pdf}
  \includegraphics[width=0.15\textwidth]{\imagepath/cell-120.pdf}
  \includegraphics[width=0.15\textwidth]{\imagepath/cell-600.pdf}
  \caption{The $5$-, $8$-, $16$-, $24$-, $120$- and $600$-cell
    tessellations of $\bbS^3$, stereographically projected to
    $\bbR^3$.}
\end{figure}

\typeout{== section2 ============================================}
\section{CMC polygons via the DPW method}
\label{sec:polygon}

\subsection{The generalized Weierstrass representation (DPW)}

Define the following loop groups (see~\cite{Pressley_Segal_1986} for
details):

\begin{itemize}
\item[]
  $\Loop$ = smooth maps (loops) from $\bbS^1_\lambda$ to
  $\matSL{2}{\bbC}$
\item[]
  $\Loop_u$ = the subgroup loops in $\Loop$ which are in $\matSU{2}{}$
  on $\bbS^1$
\item[]
  $\Loop_+$ = the subgroup loops in $\Loop$ which extend to the
  interior unit disk in $\CPone_\lambda$
\item[]
  $\Loop_{-}$ = the subgroup of loops in $\Loop$ which extend to the
  exterior unit disk in $\CPone_\lambda$
\item[]
  $\mathring\Loop_{-}$ = the subgroup of $g\in\Loop_{\NEG 1}$ such
  that $g(0)$ is upper triangular with diagonal in $\bbR^+$.
  \end{itemize}

These can be generalized to loops on a circle of radius $r\in(0,\,1)$;
see~\cite{McIntosh_1994}.

A \emph{DPW potential} $\xi$ on a Riemann surface $\Sigma$ is a
$\Loop\matsl{2}{\bbC}$ valued holomorphic differential form on
$\Sigma$ with $\xi = \sum_{k=\NEG 1}^\infty \xi_k\lambda^k$, $\det
\xi_{\NEG 1}=0$. A meromorphic DPW potential is defined analogously.

A CMC surface is constructed from a DPW potential as follows.  Let
$\Phi$ be the \emph{holomorphic frame} solving $\deriv \Phi = \Phi
\xi$; $\Phi$ generally has monodromy.  Let $\Phi = FB \in
\Loop_u\mathring\Loop_+$ be the Iwasawa factorization into the
\emph{unitary frame} $F$ and positive part $B$
(see~\cite{Pressley_Segal_1986}); for our case of $\matSU{2}{}$ this
factorization always exists.  The CMC surface is constructed via
the formulas first obtained in \cite{Bobenko_1991:cmcsurfaces}:
\begin{subequations}\label{eq:extrinsic_closing_unitary_potential}
\begin{align}
  \bbS^3 &: \quad f(\lambda_0,\,\lambda_1) = F(\lambda_0) F^{\NEG
    1}(\lambda_1) \spacecomma\quad \lambda_0,\,\lambda_1\in\bbS^1
  \\ \bbH^3 &: \quad f(\lambda_0,\,\lambda_1) = F(\lambda_0) F^{\NEG
    1}(\lambda_1) \spacecomma\quad \lambda_1 = \lambda_0^{\NEG
    1}\in\bbC\setminus\bbS^1 \\ \bbR^3 &: \quad f(\lambda_0) =
  \NEG\tfrac{2}{H}\dot F(\lambda_0) F(\lambda_0)^{\NEG 1}
  \spacecomma\quad \lambda_0 = 1 \text{ or } \NEG 1
\end{align}
\end{subequations}
where in the case of $\bbR^3$ the dot denotes the derivative with
respect to $\theta$, $\lambda = e^{\imi\theta}$. The unitary frame $F$
yields a unitary potential $\mu=F^{-1}\deriv F$ which is well-defined
on the (Riemann) surface as opposed to $F$ which is well-defined only
on the universal covering.  The unitary potential is also known as the
associated family of flat connections, see~\cite{Heller_2013} and the
literature therein.

\vspace{0.0975cm}\begin{textblock*}{16cm}(14.5cm,13.8cm)
\begin{minipage}{5cm}
\begin{figure}
  \ximage[1]
    {r3-torus-3-bulge.pdf}
    {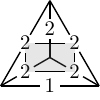}
    {\scriptsize (\textsc{a}) torus 3 (bulge) \vspace{0.35cm}}
  \ximage[1]
    {r3-torus-3-neck.pdf}
    {diagram/tet-i-r3-torus-3.pdf}
    {\scriptsize (\textsc{b}) torus 3 (neck) \vspace{0.35cm}}
  \captionsetup{margin=0cm}
  \caption{\label{fig:r3-torus-a}CMC tori in $\bbR^3$ with Delaunay ends}
\end{figure}
\end{minipage}
\end{textblock*}

A DPW potential $\xi = \sum_{k=\NEG 1}^\infty \xi_k\lambda^k$ is
\emph{adapted} if $\xi_{\NEG 1}$ is upper triangular (and hence has
zero diagonal because $\det \xi_{\NEG 1}=0$).  For adapted DPW
potentials, the \emph{Hopf differential} is $Q = \DOT{\xi_{\NEG
    1}}{\xi_0}$.  For non-adapted potentials, the formula for $Q$ is
more complicated.

If $\xi$ is holomorphic at $z_0$, the induced CMC surface is
\emph{immersed} at $z_0$ if and only if $\xi_{\NEG 1}$ does not vanish
at $z_0$.

This representation differs from the original
representation~\cite{Dorfmeister_Pedit_Wu_1998} in that the potential
and loops are not twisted.  It is slightly looser than the
representation~\cite{Schmitt_Kilian_Kobayashi_Rossman_2007} in that it
does not require the DPW potential to be adapted.

\subsection{Delaunay ends}
\label{sec:delaunay}

A \emph{Delaunay eigenvalue} is
\begin{equation}
  \label{eq:delaunay-eigenvalue}
  \nu = \half \sqrt{ 1 + \lambda^{\NEG 1}(\lambda - \lambda_0)(\lambda
    - \lambda_1)w} \spaceperiod
  \end{equation}
where the evaluation points $\lambda_0$, $\lambda_1$ and the \emph{end
weight} $w\in\Rstar$ chosen so that $\nu$ is real on $\bbS^1$.  A DPW
potential with a simple pole, unitary monodromy, and Delaunay
eigenvalues of the residue induces a surface asymptotic to a half
Delaunay cylinder~\cite{Kilian_Rossman_Schmitt_2008}.

To construct surfaces, two types of closing conditions must be
satisfied by a DPW potential:
\begin{itemize}
  \item
  The \emph{intrinsic closing condition} is the condition that the
  monodromy group is unitarizable on $\bbS^1$ (or more generally,
  $r$-unitarizable on an circle of radius $r\in(0,\,1)$).  For the
  Fuchsian DPW potentials, this condition is not directly satisfiable
  except in the case of $3$ or $2$ geometric poles; more than $3$
  requires the \emph{unitary flow}, which by definition preserves the
  intrinsic closing condition.
  \item
  The \emph{extrinsic closing conditions} are conditions on the DPW
  potential on the monodromy at the evaluation points, chosen to
  control the desired geometry of the surface via
  \eqref{eq:extrinsic_closing_unitary_potential}.  For surfaces
  constructed via tessellations these conditions are given
  in~\cref{thm:cmc-polygon}.
\end{itemize}

\subsection{Gauge}
\label{sec:gauge}

Consider a holomorphic DPW potential $\xi$ and a holomorphic map
$g\colon\Sigma\to\Loop$.  The gauge action is
\begin{equation}
\xi \mapsto g^{\NEG 1}\xi g + g^{\NEG 1}\deriv g \spaceperiod
\end{equation}
The point of the gauge action is that if $\deriv\Phi = \Phi \xi$ then
$\deriv(\Phi g) = (\Phi g) (\gauge{\xi}{g})$.  We allow gauges to have
monodromy $\pm\id$ along paths; such multivalued gauges nevertheless
map single-values potentials to single-valued potentials.

A $\emph{DPW gauge}$ is one which maps DPW potentials to DPW
potentials, that is, $g\colon\Sigma\to\Loop_+$ is holomorphic in
$\lambda$.  If $\xi$ is a DPW potential and $g$ is DPW gauge, then
$\xi$ and $\gauge{\xi}{g}$ induce the same surface in the sense that
$\Phi$ and $\Phi g$ do.  A DPW gauge $g$ is adapted if it preserves
adapted DPW potentials, that is, $\left.g\right|_{\lambda=0}$ is upper
triangular.

Let $\xi$ be a holomorphic DPW potential on $\Sigma.$ A meromorphic
DPW gauge is given by a meromorphic map $g\colon \Sigma\to\Lambda_+.$
Then, $\gauge{\xi}{g}$ is a meromorphic DPW potential with so-called
apparent singularities at the singular points of $g$.  In general, the
singularities of a meromorphic DPW potential are not apparent.

\subsection{Spin}
\label{sec:spin}
For a DPW gauge $g$ define the group homomorphism
\begin{subequations}
\begin{gather}
  \spin: H_1(\Sigma)\to\bbZ_2 = \{\pm 1\} \spacecomma \\ \spin_\gamma
  g = \begin{cases} +1 & \text{if $g$ has monodromy $+1$ along
      $\gamma$}\\ -1 & \text{if $g$ has monodromy $-1$ along $\gamma$}
\end{cases}
\end{gather}
\end{subequations}
that is, $\spin_\gamma g = +1$ (resp.~$\NEG 1$) if $g$ returns to
itself (resp.~its negative) along $\gamma$.  Then $\spin_\gamma gh =
\spin_\gamma g \cdot \spin_\gamma h$.

To define $\spin$ for DPW potentials consider the double cover
\begin{equation}
  \label{eq:double-cover}
\bbC^2\setminus\{0\}\to\{x\in\matsl{2}{\bbC}\setminus\{0\}\suchthat
\det x = 0\} \spacecomma\quad
\begin{smatrix}u\\v\end{smatrix}
  \mapsto
\begin{smatrix}u \\ v\end{smatrix}
  \begin{smatrix}\NEG v & u\end{smatrix}
    \spaceperiod
\end{equation}
For a DPW potential $\xi$ on a Riemann surface $\Sigma$, let
$\xi_{\NEG 1}$ be its $\lambda^{\NEG 1}$ coefficient.  Define the
group homomorphism
\begin{subequations}
\begin{gather}
  \spin: H_1(\Sigma)\to\bbZ_2 = \{\pm 1\} \spacecomma \\ \spin_\gamma
  \xi =
\begin{cases}
  +1 & \text{if the lift of $\xi_{\NEG 1}$ along $\gamma$ is a closed
    cuvve}\\ -1 & \text{otherwise}
 \end{cases}
\end{gather}
\end{subequations}

That is, $\spin_\gamma \xi$ is $+1$ (resp.~$\NEG 1$) if the lift of
$\xi_{\NEG 1}$ returns to itself (resp.~its negative) along $\gamma$.
Then $\spin_\gamma \gauge{\xi}{g} = \spin_\gamma \xi \cdot
\spin_\gamma g$.

\vspace{0.0975cm}\begin{textblock*}{16cm}(14.5cm,7.5cm)
\begin{minipage}{5cm}
\begin{figure}
  \ximage[1]
    {r3-torus-4-bulge.pdf}
    {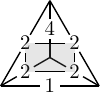}
    {\scriptsize (\textsc{a}) torus 4 (bulge) \vspace{0.35cm}}
  \ximage[1]
    {r3-torus-6-bulge.pdf}
    {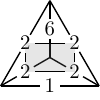}
    {\scriptsize (\textsc{b}) torus 6 (bulge) \vspace{0.35cm}}
  \ximage[1]
    {r3-torus-8-bulge.pdf}
    {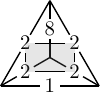}
    {\scriptsize (\textsc{c}) torus 8 (bulge) \vspace{0.35cm}}
  \captionsetup{margin=0cm}
  \caption{\label{fig:r3-torus-b}CMC tori in $\bbR^3$ with Delaunay ends}
\end{figure}
\end{minipage}
\end{textblock*}

The spin can similarly be defined for unitary potentials using the
lift of the coefficient of $\lambda^{\NEG 1}$.  If $\Phi = FB$ is the
holomorphic frame with potential $\xi$, and $F$ is the corresponding
unitary frame with potential $\eta$, then $\spin\eta = \spin\xi$
because $\spin B|_{\lambda=0}= 1$ since
$B|_{\lambda=0}\in\mathring\Lambda_+$.

A geometric interpretation for the spin of a potential can be given in
terms of a coordinate frame, that is, a unitary frame $G$ satisfying
\begin{equation}
N = G e_0 G^{\NEG 1} \spacecomma\quad f_x/v = G e_1 G^{\NEG 1}
\spacecomma\quad f_y/v = G e_2 G^{\NEG 1}
\end{equation}
where $e_0,\,e_1,\,e_2$ is a positively oriented orthonormal basis for
$\matsu{2}{}$, $f$ is the CMC immersion, $v$ is the metric of $f$, and
$N$ is its normal.  Then $\spin \xi = \spin u$, where $u = F^{\NEG
  1}G$ is the gauge between the unitary frame $F$ and a coordinate
frame $G$.

Consider a meromorphic DPW potential $\xi$ on $\Sigma$.  For
$z\in\Sigma$ write $\spin_z\xi$ to mean $\spin_\gamma\xi$ along a
small circle $\gamma$ encircling $z$. If $\xi$ is regular at $z$ then
$\spin_z\xi=1.$ For a DPW potential $\xi$ on $\CPone$ with finitely
many singularities we have the total spin
\begin{equation}
\label{eq:total-spin}
\prod_{z\in\CPone}\spin_z\xi = 1 \spaceperiod
\end{equation}

As an application of spin, when we construct CMC polygons whose
boundaries reflect in planes in \cref{sec:cmc-polygon}, the spin is
used to distinguish the internal and external dihedral angles of the
planes.

\subsection{Symmetry}
\label{sec:symmetry}

The following theorem and lemma detail how a symmetry of the potential
descends to a symmetry of the meromorphic frame, the unitary frame,
and the CMC immersion via
\eqref{eq:extrinsic_closing_unitary_potential}.

\begin{theorem}
  \label{thm:symmetry}
  \theoremname{Symmetry} Let $\xi$ be a DPW potential.
  \begin{enumerate}
  \item
    \label{item:symmetry-1}
    If for a holomorphic automorphism $\tau$ of the domain,
    $\tau^\ast\xi = \gauge{\xi}{g}$, then $\tau^\ast\Phi = R \Phi g$
    for some $R\in\Loop$.  If $R$ is unitary, then the CMC immersion
    has the orientation preserving symmetry
    \begin{subequations}
      \label{eq:symmetry-1}
      \begin{align}
      \text{$\bbS^3$ and $\bbH^3$} &:\quad \tau^\ast f = R f R^{\NEG
        1} \\ \text{$\bbR^3$} &:\quad \tau^\ast f = R f R^{\NEG 1} -
      \tfrac{2}{H} \dot{R}{R}\spacecomma
    \end{align}
    \end{subequations}
    where in the case of $\bbR^3$ the dot denotes the derivative with
    respect to $\theta$, $\lambda = e^{\imi\theta}$.
  \item
    \label{item:symmetry-2}
    If for an antiholomorphic automorphism $\tau$ of the domain,
    $\tau^\ast\xi(\lambda) = \gauge{\ol{\xi(\ol{\lambda})}}{g}$, then
    $\tau^\ast\Phi(\lambda) = R \ol{\Phi(\ol{\lambda})} g$ for some
    $R\in\Loop$.  If $R$ is unitary, then the CMC immersion has the
    orientation reversing symmetry
    \begin{subequations}
      \label{eq:symmetry-2}
    \begin{align}
      \text{$\bbS^3$ and $\bbH^3$}&:\quad \tau^\ast
      f(\lambda_0,\,\lambda_1) = R
      \ol{f(\ol{\lambda_0},\,\ol{\lambda_1})} R^{\NEG 1}\,,
      \\ \text{$\bbR^3$} &:\quad \tau^\ast f(\lambda_0) = -R
      \ol{f(\ol{\lambda_0})} R^{\NEG 1} - \tfrac{2}{H} \dot{R}{R}
      \spaceperiod
    \end{align}
    \end{subequations}
    \end{enumerate}
\end{theorem}

In the orientation reversing case of the above theorem, the
symmetry~\eqref{eq:symmetry-2} relates two associate CMC surfaces,
which are the same surface if $\lambda_1 = \ol{\lambda_0}$ (for
$\bbS^3$), $\lambda_0\in\bbR$ and $\lambda_1\in\bbR$ (for $\bbH^3$)
and $\lambda_0\in\{\pm 1\}$ (for $\bbR^3$).

\Cref{thm:symmetry} is of limited use without the knowledge that $R$
in that theorem is unitary.  One necessary condition that $R$ is
unitary is given in the following lemma:

\begin{lemma}
  \label{lem:symmetry}
  \theoremname{Unitarity of $R$} If $\xi$ in~\cref{thm:symmetry}
  extends to $\bbS^1_\lambda$ and has irreducible unitary monodromy,
  then $R$ in that theorem is unitary.
\end{lemma}
\begin{proof}
  Let $f:U\to\Sigma$ be the universal cover, and $\tau$ a lift of
  $\tau$ to the universal cover, so $f\tau = \tau\circ f$.  Let
  $\sigma$ be a deck transformation, so $f\sigma = f$.  Then
  $f\tau\sigma\tau^{\NEG 1} = f$ implying that $\tau\sigma\tau^{\NEG
    1}$ is a deck transformation.

  Let $M_\sigma$ the monodromy of $\Phi$ with respect to $\sigma$.
  and $\spin_{\sigma}\id\in\{\pm\id\}$ the monodromy of $g$ with
  respect to $\sigma$.  For the orientation reversing case,
  \begin{equation}
    \sigma^\ast\tau^\ast \Phi = \sigma^\ast(R \ol{\Phi} g) =
    (\spin_\sigma g)R \ol{M}_{\sigma} \ol{\Phi} g = (\spin_\sigma g)R
    \ol{M}_{\sigma} R^{\NEG 1} \tau^\ast \Phi
  \end{equation}
  so
  \begin{equation}
    {\tau^{\NEG 1}}^\ast\sigma^\ast\tau^\ast \Phi = (\spin_\sigma g)R
    \ol{M}_{\sigma} R^{\NEG 1} \Phi \spaceperiod
  \end{equation}
  Since $\tau\sigma\tau^{\NEG 1}$ is a deck transformation, its
  monodromy is given by $N_\sigma \coloneq (\spin_\sigma g)R
  \ol{M}_{\sigma} R^{\NEG 1}$.  Since by assumption the monodromy
  group is irreducible and unitary, then $N_\sigma\in\Loop_u$ for
  every deck transformation $\sigma$.  Using that $\xi$ extends to
  $\bbS^1_\lambda$, this implies $R\in\Loop_u$.

  The proof for the orientation preserving case is the same without
  the overline.
\end{proof}

\subsection{CMC polygons}
\label{sec:cmc-polygon}

Let $\bbR\cup\{\infty\}$ be divided into $n$ segments $s_1,\dots,s_n$
at $n$ distinct consecutive points $z_{ij}$ dividing $s_i$ and $s_j$.
Let $\xi$ be a meromorphic DPW potential on $\CPone$ with
singularities at these points $z_{ij}$.  With $b$ a basepoint in the
upper halfplane, for $i,\,j\in\{1,\dots,n\}$ let $\gamma_{ij}$ be a
simple closed counterclockwise curve based at $b$ which crosses the
segments $s_i$ and $s_j$, and let $M_{ij}$, $i,\,j\in\{1,\dots,n\}$,
$i < j$ be the monodromy along $\gamma_{ij}$.  The $n$ \emph{local}
monodromies are those along paths which enclose one singularity; the
remaining monodromies are called \emph{global}.

\vspace{0.0975cm}\begin{textblock*}{16cm}(14.5cm,7.5cm)
\begin{minipage}{5cm}
\begin{figure}
  \ximage[1]
    {r3-torus-4-neck.pdf}
    {diagram/tet-i-r3-torus-4.pdf}
    {\scriptsize (\textsc{a}) torus 4 (neck) \vspace{0.35cm}}
  \ximage[1]
    {r3-torus-6-neck.pdf}
    {diagram/tet-i-r3-torus-6.pdf}
    {\scriptsize (\textsc{b}) torus 6 (neck) \vspace{0.35cm}}
  \ximage[1]
    {r3-torus-8-neck.pdf}
    {diagram/tet-i-r3-torus-8.pdf}
    {\scriptsize (\textsc{c}) torus 8 (neck) \vspace{0.35cm}}
  \captionsetup{margin=0cm}
  \caption{\label{fig:r3-torus-c}CMC tori in $\bbR^3$ with Delaunay ends}
\end{figure}
\end{minipage}
\end{textblock*}

\begin{theorem}
  \theoremname{CMC polygons}
  \label{thm:cmc-polygon}
  Let $\xi$ be a meromorphic DPW potential satisfying the conditions
  of \cref{lem:symmetry} with $n$ singularities on
  $\bbR\cup\{\infty\}$ as above. Assume $\xi$ admits the reflection
  symmetry $\tau^\ast\xi(\lambda) = \ol{\xi(\ol{\lambda})}$ for
  $\tau(z) = \ol{z}$.  Let $\theta_{ij}\in[0,\,\pi]$,
  $i,\,j\in\{1,\dots,n\}$, $i<j$.  If the monodromies $M_{ij}$ satisfy
  \begin{equation}
  \half \tr M_{ij}|_{\lambda_0} = \NEG (\spin_{\gamma_{ij}}\xi)
  \cos\theta_{ij} \spacecomma\quad
  i,\,j\in\{1,\dots,n\}\spacecomma\quad i < j
  \end{equation}
  then the CMC surface induced by $\xi$ with the upper halfplane as
  domain is a $n$-gon whose boundaries reflect in $n$ planes
  (respectively totally geodesic spheres) $P_1,\dots,P_n$, with
  internal dihedral angles $\theta_{ij}$ between $P_i$ and $P_j$.
\end{theorem}

\begin{proof}
  Let $F$ be the unitary frame, $G$ a coordinate frame, with respect
  to a basis $\hat e_0$, $\hat e_1$, $\hat e_2\in\matsl{2}{\bbC}$ and
  $u$ the unitary $\lambda$-independent gauge between them, so $F =
  Gu$.  By the proof of \cref{thm:symmetry}\eqref{item:symmetry-2},
  $\tau_k^\ast F = P_k \ol{F}$, $k\in\{1,\dots,n\}$ so
  \begin{equation}
  \tau_k^\ast G = P_k \ol{G} Q_k^{\NEG 1} \spacecomma\quad \tau^\ast u
  = Q_k \ol{u} \spaceperiod
  \end{equation}
  Then with $\rho_{ij}$ and $\sigma_{ij} = \spin_{\gamma_{ij}}u$,
  \begin{subequations}
  \begin{align}
    u &= \tau_k^\ast\tau_k^\ast u = Q_k\ol{Q}_{k}u \quad\implies\quad
    Q_k\ol{Q}_k = \id\\ \sigma_{ij}u &= \rho^\ast u = \tau_j\tau_iu =
    Q_jQ_i^{\NEG 1}u \quad\implies\quad Q_{j} = \sigma_{ij}Q_{i}
    \spaceperiod
  \end{align}
  \end{subequations}
  With $p_k$ a fixed point of $\tau_k$, define $u_k = u(p_k)$ so $u_k
  = Q_k \ol{u}_k$.  Then for $i,\,j\in\{1,\dots,n\}$, $i\ne j$, $u_j =
  Q_j \ol{u_j} = \sigma_{ij}Q_i \ol{u_j}$.  Thus $u_i^{\NEG 1}u_j =
  \sigma_{ij}\ol{u_i^{\NEG 1}u_j}$.  so
  \begin{subequations}
  \begin{align}
    \sigma_{ij} = 1 &:\quad u_i^{\NEG 1}u_j = \ol{u_i^{\NEG 1}u_j}
    \quad\implies u_i^{\NEG 1}u_j e_1 = e_1 u_i^{\NEG
      1}u_j\\ \sigma_{ij} = \NEG 1 &:\quad u_i^{\NEG 1}u_j =
    -\ol{u_i^{\NEG 1}u_j} \quad\implies u_i^{\NEG 1}u_j e_1 = -e_1
    u_i^{\NEG 1}u_j \spaceperiod
  \end{align}
  \end{subequations}
  Hence $u_i^{\NEG 1}u_j = \sigma_{ij} e_1 u_i^{\NEG 1}u_j$ so $u_j
  e_1 u_j^{\NEG 1} = \sigma_{ij} u_i e_1 u_i^{\NEG 1}$.

  Let $e_1 = \begin{smatrix}0 & 1\\\NEG 1 & 0\end{smatrix}$.  Since
    $G$ is a coordinate frame, we have ${(f_y)}_k = G_k \hat e_2
    G_k^{\NEG 1}$.  Define $\sigma_k\in\{\pm 1\}$ by $P_k e_1 =
    \sigma_k {(f_y)}_k$.  Then
  \begin{equation}
  N_k \coloneq P_k e_1 = F_k \ol{F}_k^{\NEG 1}e_1 = F_k e_1 F_k^{\NEG
    1} = G_k u_k e_1 u_k^{\NEG 1} G_k^{\NEG 1} = \sigma_k G_k \hat e_2
  G_k^{\NEG 1}
  \end{equation}
  so
  \begin{equation}
  u_k e_1 u_k^{\NEG 1} = \sigma_k \hat e_2 \spaceperiod
  \end{equation}
  Hence $\sigma_i\sigma_j = \sigma_{ij}$.

  Since $\tfrac{\deriv}{\deriv y}$ is pointing into the upper half
  plane, ${(f_y)}_k$ is an internal normal to the plane.  Thus
  $\sigma_k = 1$ if $N_k$ is internal, and $\sigma_k = \NEG 1$ if
  $N_k$ is external.  Thus $\sigma_{ij}=\sigma_i\sigma_j = 1$ if and
  only if $N_i$ and $N_j$ are both internal or both internal, and
  $\sigma_{ij}=\sigma_i\sigma_j = \NEG 1$ if and only if one of $N_i$
  and $N_j$ is internal and one external.

  This means
  \begin{equation}\DOT{N_i}{N_j} = -\sigma_{ij} \cos\theta_{ij}
  \end{equation}
  where $\theta_{ij}$ is the internal angle between planes $i$ and
  $j$.  Since $M_{ij} = P_j P_i^{\NEG 1}$, then
  \begin{equation}
  \half\trace M_{ij}|_{\lambda=\lambda_0} = \DOT{N_i}{N_j} =
  -\sigma_{ij} \cos\theta_{ij} \spaceperiod \qedhere
  \end{equation}
\end{proof}

\begin{remark}
  \theoremname{Translational monodromy} Planes with dihedral angle
  $\theta_{ij}=\pi$ are parallel.  Constraining the planes to coincide
  (for example, for the CMC torus with ends in $\bbR^3$) requires that
  the extrinsic conditions of \cref{thm:cmc-polygon} be augmented with
  the additional condition
  \begin{equation}
  \tfrac{\deriv}{\deriv\lambda} M_{ij}|_{\lambda=\lambda_0} = 0
  \spaceperiod
  \end{equation}
\end{remark}

It remains to control the vertices of the CMC polygon constructed in
\cref{thm:cmc-polygon}.  For this we use a DPW potential with simple
poles:

\begin{theorem}
  \label{thm:cmc-polygon-vertex}
  \theoremname{CMC polygon vertices} Let $\xi$ a DPW potential as in
  \cref{thm:cmc-polygon}, $z_k$ a simple pole of $\xi$ on $\bbR$, and
  $\nu$ the eigenvalue of $\res_{z=z_k}\xi$.
  \begin{enumerate}
  \item
    \label{item:cmc-polygon-vertex-1}
    If $\nu=1/(2n)$ or $\nu=\half-1/(2n)$, $n\in\bbN_{\ge 2}$, and
    $\xi_{-1}$ has a simple pole at $z_k$, then the CMC surface
    constructed from $\xi$ with $2n$ reflections around the vertex is
    immersed at the vertex.
  \item
    \label{item:cmc-polygon-vertex-2}
    If $\nu = \nu_{\mathrm{Del}}/n$, $n\in\bbN_{\ge 2}$, where
    $\nu_{\mathrm{Del}}$ is a Delaunay
    eigenvalue~\eqref{eq:delaunay-eigenvalue}, then the CMC surface
    constructed from $\xi$ with $2n$ reflections around the vertex is
    a once-wrapped Delaunay end.
    \end{enumerate}
\end{theorem}

\vspace{0.0975cm}\begin{textblock*}{16cm}(14.5cm,7.5cm)
\begin{minipage}{5cm}
\begin{figure}
  \ximage[1]
    {r3-platonic-tet-bulge.pdf}
    {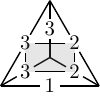}
    {\scriptsize (\textsc{a}) tetrahedron (bulge) \vspace{0.35cm}}
  \ximage[1]
    {r3-platonic-tet-neck.pdf}
    {diagram/tet-i-r3-platonic-tet.pdf}
    {\scriptsize (\textsc{b}) tetrahedron (neck) \vspace{0.35cm}}
  \ximage[1]
    {r3-platonic-oct-bulge.pdf}
    {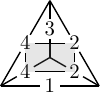}
    {\scriptsize (\textsc{c}) octahedron \vspace{0.35cm}}
  \captionsetup{margin=0cm}
  \caption{\label{fig:r3-platonic-a}CMC surfaces in $\bbR^3$ with Platonic symmetry and Delaunay ends}
\end{figure}
\end{minipage}
\end{textblock*}

\begin{proof}
  Assuming $z_k=0$, write $\xi = A_{\NEG 1}\deriv z/z + A_0\deriv z +
  \dots$.  The pullback with respect to the local covering map $f(w) =
  w^n$ is
  \begin{equation}
  f^\ast \xi = nA_{\NEG 1} \deriv w/w + nA_0 w^{n-1}\deriv w + \dots
  \spaceperiod
  \end{equation}

  Proof of~\eqref{item:cmc-polygon-vertex-1}: Since by assumption
  $A_{-1}$ has a pole at $\lambda = 0$, by a $z$-independent local
  gauge of $\xi$ it may be assumed
  \begin{equation}
  A_{\NEG 1} = \begin{bmatrix}1/(2n) & \lambda^{\NEG 1} \\ 0 & \NEG
    1/(2n)\end{bmatrix} \spaceperiod
  \end{equation}
  Then the local gauge $g=\diag(w^{\NEG 1/2},\,w^{1/2})$ removes the
  simple pole of $f^\ast \xi$ at $w=0$, and
  \begin{equation}
  \gauge{(f^\ast\xi)}{g} =
  \begin{bmatrix}0 & n\lambda^{-1}\\0 & 1\end{bmatrix}\deriv w + \dots
  \spaceperiod
  \end{equation}
  Since $\gauge{(f^\ast \xi)}{g}$ is holomorphic at $w=0$ and its
  $\lambda^{-1}$ coefficient does not vanish at $w=0$, then the CMC
  surface induced by $\gauge{(f^\ast\xi)}{g}$ is immersed at
  $w=0$. The proof for $\nu=\half-1/(2n)$ is analogous.

  Proof of~\eqref{item:cmc-polygon-vertex-2}: Since the eigenvalue of
  $A_{\NEG 1}$ is $\nu_{\mathrm{Del}}/n$, then the eigenvalue of
  $nA_{\NEG 1}$ is $\nu_{\mathrm{Del}}$, Unitary monodromy implies
  this is a once-wrapped Delaunay end.
\end{proof}

\typeout{== section3 ============================================}\section{Symmetric CMC surfaces with genus}
\label{sec:experimental-part}

\subsection{The potential}
\label{sec:potential}
By applying a M\"obius transformation we assume that the singular
points of the CMC polygon are on the unit circle.  As the fundamental
piece is a CMC quadrilateral, we restrict to the 4-punctured sphere in
the following.  We will see in \cref{sec:FuchsianDPW4} that, at least
for surfaces without Delaunay ends, we can restrict without loss of
generality to a Fuchsian DPW potential of the 4-punctured sphere. The
means it has four simple poles and no pole at $z=\infty$, and is of
the form
\begin{equation}
  \label{eq:potential}
\xi = \sum_{k=0}^3 \frac{A_k}{z-z_k}\deriv z
\end{equation}
as follows:
\begin{itemize}
\item
  The poles are $z_0\in\bbS^1$ in the open first quadrant, and
  $(z_1,\,z_2,\,z_3)$ is a permutation of $(1/z_0,\,-z_0,\,-1/z_0)$.
\item
  The residues are
  \begin{subequations}\label{symmetric-potential}
  \begin{gather}
    A_0 = \begin{bmatrix}y &
      \lambda^{-1}p\\\frac{\lambda(\nu_0^2-y^2)}{p} & -y\end{bmatrix}
      \spacecomma\quad A_2 = \begin{bmatrix}-y & \frac{(\nu_1^2 -
          y^2)}{x}\\ x & y\end{bmatrix} \spacecomma\ \\ A_1 = \sigma
        A_0 \sigma^{-1} \spacecomma\quad A_3 = \sigma A_2
        \sigma^{-1}\spacecomma\quad \sigma = \diag(\imi,\,-\imi)
        \spaceperiod
  \end{gather}
  \end{subequations}
\item
  For surfaces without Delaunay ends, the eigenvalue $\nu_0$ of $A_0$
  and $A_1$ and $\nu_1$ of $A_2$ and $A_3$ are constants in
  $(0,\,1/2)$.  For surfaces with Delaunay ends, $\nu_0$ is of the
  form $c\nu_{\mathrm{Del}}$ with $c\in(0,\,1)$ and
  $\nu_{\mathrm{Del}}$ is the eigenvalue of a Delaunay unduloid
\begin{equation}
  \nu_{\mathrm{Del}} = \half\sqrt{1 + \fourth\lambda^{\NEG
      1}{(\lambda-1)}^2 w } \spacecomma\quad w\in (0,\,1] \spaceperiod
\end{equation}
\item
  The accessory parameters $x$ and $y$ are holomorphic functions of
  $\lambda$ on an open disk $\calD_r$ of radius $r>0$ centered at the
  origin satisfying $x(\ol{\lambda}) = \ol{x(\lambda)}$ and
  $y(\ol{\lambda}) = \ol{y(\lambda)}$.  The function $p$ is a monic
  polynomial in $\lambda$ satisfying $p(\ol{\lambda}) =
  \ol{p(\lambda)}$.
\item
  The quotients $\lambda(\nu_0^2 - y^2)/p$ and $(\nu_1^2 - y^2)/x$ are
  holomorphic functions of $\lambda$ on $\calD_r$.
\end{itemize}

The need for the last condition is as follows.  The unitary flow,
which preserves the unitarizability of the monodromy of $\xi$, is
implemented by evaluating the monodromy of $\xi$ directly on the unit
circle, and not by the numerically more problematic procedure of
computing the monodromy on an $r<1$ circle and then extending it to
the unit circle.

Let $M_k$ ($k\in\{0,\dots,3\}$) be the local monodromy around $z_k$
based at $z=1$.

The surfaces are constructed by running the unitary flow (see
\cref{sec:unitary-flow} below) so that at the end of the flow for
$k=0,1$
\begin{subequations}
\begin{gather}
\nu_0|_{\lambda=\lambda_k} = \tfrac{1}{2n_0} \spacecomma\quad
\nu_1|_{\lambda=\lambda_k} = \half - \tfrac{1}{2n_1} \spacecomma\quad
\\ \trace M_0M_1|_{\lambda=\lambda_k} = -\cos\tfrac{2\pi}{r}
\spacecomma\quad \half \trace M_1M_2|_{\lambda=\lambda_k} =
\cos\tfrac{\pi}{s}\,.
\end{gather}
\end{subequations}

\noindent
\begin{wrapfigure}{r}{0.1\linewidth}
\includegraphics[scale=1.25]{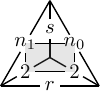}
\end{wrapfigure}
Then, by \cref{thm:cmc-polygon,thm:cmc-polygon-vertex}, the unit
$z$-disk maps to a CMC quadrilateral whose edges reflect in planes
(respectively geodesic 2-spheres) with internal dihedral angles
specified by the figure at right, and whose vertices after these
reflections are either immersed points or once-wrapped Delaunay ends.

In the special case $\nu_0 + \nu_1 = \half$, the surfaces are
constructed by running the unitary flow so that at the end of the flow
for $k=0,1$
\begin{subequations}
\begin{gather}
\nu_0|_{\lambda=\lambda_k} = \tfrac{1}{2n_0} \spacecomma\quad
\nu_1|_{\lambda=\lambda_k} = \half - \tfrac{1}{2n_0} \spacecomma\quad
\\ \trace M_0M_1|_{\lambda=\lambda_k} = -\cos\tfrac{2\pi}{r}
\spacecomma\quad \half \trace M_1M_2|_{\lambda=\lambda_k} =
\cos\tfrac{2\pi}{s}
\end{gather}
\end{subequations}

\noindent
\begin{wrapfigure}{r}{0.1\linewidth}
\includegraphics[scale=1.25]{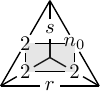}
\end{wrapfigure}
Then the quarter disk in the first quadrant maps to a CMC
quadrilateral whose edges reflect in planes with internal dihedral
angles specified by the figure at right, and whose vertices after
these reflections are immersed.

\subsection{The initial condition}
\label{sec:initial}

\subsubsection{The initial condition}
The initial condition for the unitary flow is a potential $\xi_0$ of
the form in \cref{sec:potential} with eigenvalues
$\nu_0=\nu_1=\tfrac{1}{4}$ with unitary monodromy on $\bbS^1$ which
induces a Delaunay surface.

\vspace{0.0975cm}\begin{textblock*}{16cm}(14.5cm,7.5cm)
\begin{minipage}{5cm}
\begin{figure}
  \ximage[1]
    {r3-platonic-cube-bulge.pdf}
    {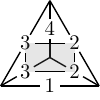}
    {\scriptsize (\textsc{a}) cube \vspace{0.35cm}}
  \ximage[1]
    {r3-platonic-ico-bulge.pdf}
    {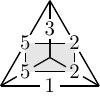}
    {\scriptsize (\textsc{b}) icosahedron \vspace{0.35cm}}
  \ximage[1]
    {r3-platonic-dodec-bulge.pdf}
    {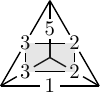}
    {\scriptsize (\textsc{c}) dodecahedron \vspace{0.35cm}}
  \captionsetup{margin=0cm}
  \caption{\label{fig:r3-platonic-b}CMC surfaces in $\bbR^3$ with Platonic symmetry and Delaunay ends}
\end{figure}
\end{minipage}
\end{textblock*}

\begin{lemma}\label{lem:tori-sd}
  \mbox{}

  \begin{enumerate}
  \item
For CMC tori of spectral genus $0$ the spectral curve
$\pi:\CPone\to\CPone$ can be chosen to be $\pi(\xi) = \xi^2$.  The
involutions are
\begin{equation}
\sigma(\xi) = -\xi \spacecomma\quad \rho(\xi) = \ol{\xi}^{-1}
\spacecomma\quad \kappa(\xi) = \ol{\xi} \spaceperiod
\end{equation}
The monodromy eigenvalues of the vacuum are
$\exp(\pm\nu_1+2\imi\pi\bbZ)$, $\exp(\pm\nu_2+2\imi\pi\bbZ)$ where
    \begin{equation}
\nu_1(\xi) \coloneq \imi\pi \frac{\xi - \xi^{-1}}{\xi_0 - \xi_0^{-1}}
\spacecomma\quad \nu_2(\xi) \coloneq \imi\pi \frac{\xi +
  \xi^{-1}}{\xi_0 + \xi_0^{-1}}
    \end{equation}
and $\pi(\pm\xi_0),\,\pi(\pm\xi_0^{-1})\in\bbS^1$ are the evaluation
points.
  \item
For some $\ell,\,m\in\bbZ^+$ the monodromy eigenvalues of a Delaunay
cylinder are $\exp(\pm\nu_1+2\imi\pi\bbZ)$,
$\exp(\pm\nu_2+2\imi\pi\bbZ)$ where
\begin{subequations}
\label{eq:delaunay-eigenvalues}
\begin{gather}
\label{eq:delaunay-eigenvalues-1}
\nu_1(u) \coloneq \imi\pi\frac{f_1(u)-f_2(u)}{f_1(u_0)-f_2(u_0)}
\spacecomma\\ \nu_2(u) = \half \ell(f_1(u) + f_2(u)) \spacecomma\quad
\nu_2(u_0) = \imi\pi m \spacecomma
\end{gather}
\end{subequations}
where $\pi(\pm u_0)$, $\pi(\pm\ol{u_0}+\half\omega_1)\in\bbS^1$ are
the evaluation points, and $f_1,\,f_2$ are as
in~\eqref{eq:delaunay-torus-functions}.
  \end{enumerate}
\end{lemma}

On the torus $\bbC/(\bbZ + \tau\bbZ)$, Let $\wp$ be the Weierstrass
function and let $\zeta\coloneq -\int \wp$.  Let
$\{\omega_1,\,\omega_2,\,\omega_3\}=\{\half,\,\half+\tfrac{\tau}{2},\,\tfrac{\tau}{2}\}$.

Define on some torus with modulus $\tau_{\mathrm{spec}}$
\begin{equation}
\label{eq:delaunay-torus-functions}
h_1(u) \coloneq \eta_1u - \omega_1\zeta(u) \spacecomma\quad h_2(u)
\coloneq f_1(u-\half\omega_1) \spaceperiod
\end{equation}

The theta function
\begin{equation}
  \theta(x,\,\tau)\coloneq \sum_{k\in\bbZ}\exp\bigl(2\imi\pi(\half
  n^2\tau + n (x-\omega_2)\bigr)
\end{equation}
is an entire function $\bbC\to\bbC$ with simple zeros at lattice
points $\bbZ + \tau\bbZ$ and no other zeros, satisfying
\begin{equation}
  \label{eq:the-gauge-theta}
  \theta(x+1) = \theta(x) \spacecomma\quad \theta(x+\tau) =
  -\exp(-2\imi\pi x)\theta(x) \spacecomma\quad \theta(\tau-x) =
  \theta(x)
\end{equation}
for all $x\in\bbC$. Define
\begin{subequations}
\begin{align}
  \half g_0(\half x) &=
  %f_0(x) - 2f_0(\half x)
  \frac{\theta'(x)}{\theta(x)} - \frac{2\theta'(\half x)}{\theta(\half
    x)} + \imi \pi \spacecomma\\ \half g_k(\half x) &= \exp\bigl(
  2\imi\pi x \tfrac{\omega_k - \ol{\omega_k}}{\tau-\ol{\tau}}\bigr)
  \frac{\theta(x+\omega_k)\theta'(0)}{\theta(x)\theta(\omega_k)}
  \spacecomma\quad k\in\{1,\,2,\,3\} \spaceperiod
\end{align}
\end{subequations}

\vspace{0.0975cm}\begin{textblock*}{16cm}(14.5cm,1cm)
\begin{minipage}{5cm}
\begin{figure}
  \ximage[1]
    {r3-fournoid-und.pdf}
    {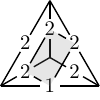}
    {\scriptsize (\textsc{a}) fournoid \vspace{0.35cm}}
  \ximage[1]
    {r3-fournoid-cyl.pdf}
    {diagram/tet-i-r3-fournoid.pdf}
    {\scriptsize (\textsc{b}) fournoid \vspace{0.35cm}}
  \ximage[1]
    {r3-fournoid.pdf}
    {diagram/tet-i-r3-fournoid.pdf}
    {\scriptsize (\textsc{c}) fournoid \vspace{0.35cm}}
  \ximage[1]
    {r3-fournoid-binoculars.pdf}
    {diagram/tet-i-r3-fournoid.pdf}
    {\scriptsize (\textsc{d}) fournoid \vspace{0.35cm}}
  \captionsetup{margin=0cm}
  \caption{\label{fig:r3-fournoid}CMC sphere in $\bbR^3$ with Delaunay ends}
\end{figure}
\end{minipage}
\end{textblock*}

The initial condition is the potential in \cref{sec:potential} with
\begin{subequations}
\begin{gather}
  \label{eq:initial-potential}
  x = \lambda(y + \nu_0)(y + \nu_1)\frac{1-u}{1+u} \spacecomma\quad
  y(b,\,a) = -\frac{\tfrac{2\imi \pi}{\tau-\ol{\tau}}(b+a) + f_0(b)}{2
    g_2(\half b)} \\ u(b) = -\frac{g_1(\half b)}{g_3(\half b)}
  \spacecomma\quad v(b,\,a) =
  \frac{\tfrac{2\imi\pi}{\tau-\ol{\tau}}(b+a) + f_0(b)}{g_2(\half b)}
  \\ a = \tfrac{2\imi\pi}{\tau-\ol{\tau}}\half h_1 \spacecomma\quad b
  = \tfrac{2\imi\pi}{\tau-\ol{\tau}}\half h_2 \\ \nu_0 = \nu_1 =
  \fourth \spacecomma\quad [z_0,\,z_1,\,z_2,\,z_3] = u(\omega_2)^2
  \spacecomma\quad p = 1 \spaceperiod
\end{gather}
\end{subequations}
The initial condition is computed numerically from
\eqref{eq:initial-potential} as Laurent series on $\bbS^1$ by
computing its Fourier coefficients. The initial data can be computed
from Lemma \eqref{lem:tori-sd} using results
from~\cite{Heller_Heller_2016}

\subsubsection{Configurations of the initial condition}
Permuting the lattice generators in the initial condition creates
different arrangements of residues of the DPW potential on the
Delaunay cylinder.  For the configurations used in this report, the
two circle arcs ($z_0,\,z_1)$ and $(z_2,\,z_3)$ are mapped to
semicircles (resp.~profile curves) on the Delaunay surface, while the
other two circle arcs ($z_1,\,z_2)$ and $(z_3,\,z_0)$ are mapped to
profile curves (resp.~semicircles) on the Delaunay surface.  The first
of these configurations is used to compute the 2-dimensional lattices
and the cubic lattices; the second is used to compute the tori and
Platonic surfaces with ends.

\subsubsection{Neck and bulge}
For the initial potential $\xi_0$ above, the poles of the Fuchsian DPW
potential are at necks of the Delaunay surface.  The initial potential
with poles at bulges is constructed as a gauge of $\xi_0$ by the gauge
$\diag( {(\lambda-\lambda_0)}^{-1/2}, {(\lambda-\lambda_0)}^{1/2} )$,
where the $\lambda_0$ is a common zero of $x$ and $y^2 - 1/16$.  This
gauge is not a DPW gauge, but a so-called dressing transformation.

\subsection{The unitary flow}
\label{sec:unitary-flow}

\subsubsection{The unitary flow}

The unitary flow is a flow through the space of potentials of
\cref{sec:potential} preserving the intrinsic closing condition. It
starts at a potential in the space with unitary monodromy, and flows
until the monodromy at the evaluation points reach some desired
extrinsic closing conditions.

Given a smooth function $F = F(t,\,\vec{x}):\bbR^{1+n}\to\bbR^n$
encoding $n$ conditions on a flow parameter $t\in[0,\,1]$ and $n$
variables $\vec{x}$, if $\det\frac{\deriv F}{\deriv x} \ne 0$, then
$x(t)$ satisfying $F(t,\,x(t))=0$ can be computed by the implicit ODE
$ \frac{\deriv F}{\deriv t} + \frac{\deriv F}{\deriv x}\frac{\deriv
  x}{\deriv t} = 0.  $ The solution for the infinite case can be
computed numerically by truncating to $F:\bbR^{1+n}\to\bbR^m$, $m\ge
n$, and solving the resulting finite dimensional ODE by least squares
methods.

The variables $\vec{x}$ parametrizing the potential consists of:
\begin{itemize}
  \item
    the conformal type $[z_0,\,z_1,\,z_2,\,z_3]$
  \item
    the local eigenvalues $\nu_0|_{\lambda=1}$ and
    $\nu_1|_{\lambda=1}$
  \item
    the end weight $w_0$
  \item
    the polynomial $p$
  \item
    the accessory parameters $x$ and $y$.\
\end{itemize}

The accessory parameters, which are holomorphic functions of
$\lambda$, are approximated by truncating their power series at
$\lambda=0$. We always assume that these functions extend to a disc
$D_r$ with $r>1$ so that the first bullet point below can be checked
directly on the unit circle.

The constraints $F=0$ are of two types:
\begin{itemize}
\item
  the intrinsic closing conditions: the halftraces $t_{ij} = \half \tr
  M_iM_j$, $i,\,j\in\{0,\dots,3\}$, $i<j$ are real on $\bbS^1$.  This
  is a necessary condition by \cref{thm:unitarizability}.
\item
  geometric constraints which choose a path through the space of
  geometric parameters to reach the desired extrinsic closing
  conditions at the end of the flow.
\end{itemize}
 By \cref{thm:unitarizability} the monodromy is unitarizable if all
 halftraces along the unit circle are real and of absolute value less
 or equal to 1.  As the components of irreducible
 $\mathrm{SL}(2,\bbC)$-representations of the 4-punctured sphere with
 real traces consists entirely of either $\mathrm{SL}(2,\bbR)$ or
 $\mathrm{SU}(2)$ representations, and since $\xi_0$ is unitarizable,
 we can ignore the condition that all traces are of absolute value
 less or equal to 1 during the unitary flow.

\vspace{0.0975cm}\begin{textblock*}{16cm}(14.5cm,7.5cm)
\begin{minipage}{5cm}
\begin{figure}
  \ximage[1]
    {s3-lawson-xi12.pdf}
    {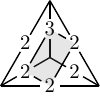}
    {\scriptsize (\textsc{a}) Lawson $\xi_{12}$ \vspace{0.35cm}}
  \ximage[1]
    {s3-lawson-xi13.pdf}
    {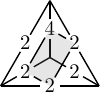}
    {\scriptsize (\textsc{b}) Lawson $\xi_{13}$ \vspace{0.35cm}}
  \ximage[1]
    {s3-lawson-xi22.pdf}
    {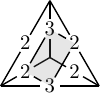}
    {\scriptsize (\textsc{c}) Lawson $\xi_{22}$ \vspace{0.35cm}}
  \captionsetup{margin=0cm}
  \caption{\label{fig:s3-lawson}Minimal Lawson surfaces in $\bbS^3$}
\end{figure}
\end{minipage}
\end{textblock*}

The intrinsic closing conditions on $\bbS^1_\lambda$ are approximated
by evaluation at finitely many equally spaced sample points on
$\bbS^1_\lambda$.  In the following we describe the other constraints
in more detail:

\subsubsection{Geometric constraints}

The simplest configuration of the geometric constraints are as
follows.  The two local and two global eigenvalues depend linearly on
the flow parameter $t$ to reach the desired values at $t=1$.  If the
surface has no ends, the end weight $w_0$ is set to $0$; otherwise it
depends linearly on $t$ starting at $0$ and reaching a heuristically
chosen value at $t=1$.  In this configuration the conformal type is
fixed during the flow.

It is possible that the flow with this simple configuration breaks
down, in which case the path must be modified in some heuristically
determined way, for example by making the conformal type depend on the
flow parameter.

In practice each geometric parameter is of one of three types:
\begin{itemize}
\item fixed during the flow
\item depending linearly on the flow parameter $t$
\item free (unconstrained).
\end{itemize}
Then the fixed variables, and the variables depending on $t$, being
computable from $t$, can be omitted from $\vec{x}$.

\subsection{Irreducibility and unitarizability}

For a subgroup $\calG\subset\matSL{2}{\bbC}$ generated by three
elements, this section proves

\begin{itemize}
\item
  a necessary and sufficient condition for the irreducibility of
  $\calG$, and
\item
  a necessary condition for the $\matSU{2}{}$ unitarizability of
  $\calG$, assuming $\calG$ is irreducible.
\end{itemize}

Here, a group $\calG$ is \emph{reducible} if all elements have a
common eigenline, and is $\matSU{2}{}$ unitarizable if there exists
$C\in\matSL{2}{\bbC}$ such that $C\calG C^{\NEG 1}\subset\matSU{2}{}$.
The methods used in the proofs can be generalized to any finitely
generated group.

The proof depends on the following \cref{lem:c3geometry}, which
determines to what extend three elements of $\bbC^3$ are determined by
their standard $\bbC^3$ inner products.

With $\DOT{-}{-}$ the standard inner product on $\bbC^3$, let $\calL =
\{v\in\bbC^3\suchthat \DOT{v}{v}=0\}$.  Let $X =
(x_0,\,x_1,\,x_2)\in\matN{3}{\bbC}$, with columns
$x_0,\,x_1,\,x_2\in\bbC^3$.  Let $W =
\transpose{X}X\in\matsym{n}{\bbC}$, so $W_{ij} = \DOT{x_i}{x_j}$.

\begin{lemma}
  \label{lem:c3geometry}
  \theoremname{Geometry of $\bbC^3$} With $X$ and $W$ as above,
  \begin{enumerate}
  \item
    \label{item:c3geometry-1}
    $\ker\transpose{X} \cap \calL = \{0\}$ if and only if $\rank W \ge
    2$.
  \item
    \label{item:c3geometry-2}
    Assuming (a), if for some $Y\in\matN{3}{\bbC}$, $\transpose{X}X =
    \transpose{Y}Y$ and $\det X = \det Y$, then there exists a unique
    $S\in\matSO{3}{\bbC}$ such that $Y = S X$.
  \end{enumerate}
\end{lemma}

\begin{proof}
  By the rank-nullity theorem applied to $\transpose{X}|_{\image X}$,
  \begin{equation}
    \label{eq:c3geometry-A}
    \rank X = \dim( \ker \transpose{X} \cap \image X ) + \rank W
  \end{equation}
  from which it follows that $\rank X \ge \rank W$, and $\rank W = 3$
  if $\rank X = 3$.

  Moreover, if $\rank X \ge 2$, then
  \begin{equation}
    \label{eq:c3geometry-B}
    \ker \transpose{X} \cap \image X = \ker \transpose{X} \cap \calL
    \spaceperiod
  \end{equation}

  To prove \eqref{item:c3geometry-1}, assume $\rank W \ge 2$, so
  $\rank X = 2$.  By~\eqref{eq:c3geometry-A}, $\dim( \ker
  \transpose{X} \cap \image X)=0$.  By~\eqref{eq:c3geometry-B}, $\ker
  \transpose{X} \cap \calL = \{0\}$.

  Conversely, assume $\ker\transpose{X}\cap\calL = \{0\}$.  Then
  $\rank X\ge 2$ because every $2$-dimensional subspace of $\bbC^3$
  intersects $\calL$.  By~\eqref{eq:c3geometry-B}, $\ker \transpose{X}
  \cap \calL = \{0\}$.  By~\eqref{eq:c3geometry-A}, $\rank W = \rank X
  \ge 2$.

  To prove \eqref{item:c3geometry-2} in the case $\rank W = 3$, since
  $\rank X = \rank Y = 3$, define $S \coloneq YX^{-1}$.  Then
  $S\in\matSO{3}{\bbC}$ by $\transpose{X}X = \transpose{Y}{Y}$ and
  $\det X = \det Y$.

  To prove \eqref{item:c3geometry-2} in the case $\rank W = 2$, let
  $x_a,\,x_b$ be two independent columns of $X$ and let $\hat X =
  (x_a,\,x_b,\,x_a\times x_b)$ and $\hat Y = (y_a,\,y_b,\,y_a\times
  y_b)$.  Since $x_a\times x_b\in\ker \transpose{X}$, then by the
  assumption and \cref{lem:c3geometry}\eqref{item:c3geometry-1},
  $x_a\times x_b\not\in\calL$.  Then $\det\hat X = \DOT{x_a\times
    x_b}{x_a\times x_b} \ne 0$ so $\rank \hat X = 3$.  Moreover, since
  $\DOT{x_a\times x_b}{x_a\times x_b} = \DOT{y_a\times y_b}{y_a\times
    y_b}$, then $\transpose{\hat{X}}\hat{X} =
  \transpose{\hat{Y}}\hat{Y}$.  Then $S \coloneq \hat{Y} \hat{X}^{-1}$
  is in $\matSO{3}{\bbC}$, and $Y = S X$.
\end{proof}

\vspace{0.0975cm}\begin{textblock*}{16cm}(14.5cm,7.5cm)
\begin{minipage}{5cm}
\begin{figure}
  \ximage[1]
    {s3-platonic-tet.pdf}
    {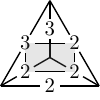}
    {\scriptsize (\textsc{a}) tetrahedron \vspace{0.35cm}}
  \ximage[1]
    {s3-platonic-oct.pdf}
    {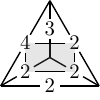}
    {\scriptsize (\textsc{b}) octahedron \vspace{0.35cm}}
  \ximage[1]
    {s3-platonic-cube.pdf}
    {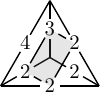}
    {\scriptsize (\textsc{c}) cube \vspace{0.35cm}}
  \captionsetup{margin=0cm}
  \caption{\label{fig:s3-platonic-a}Minimal surfaces in $\bbS^3$ with Platonic symmetry}
\end{figure}
\end{minipage}
\end{textblock*}

Identify $\bbC^4$ with $\matgl{2}{\bbC}$ by identifying the standard
basis $E_0,\,E_1,\,E_2,\,E_3$ with $\id,\,e_0,\,e_1,\,e_2$, where
\begin{equation}
e_0 \coloneq
\begin{smatrix}\imi & 0\\0 & \NEG\imi\end{smatrix}
\spacecomma\quad e_1 \coloneq
  \begin{smatrix}0 & 1\\-1 & 0\end{smatrix}
\spacecomma\quad e_2 \coloneq
  \begin{smatrix}0 & \imi\\\NEG\imi & 0\end{smatrix}
\spaceperiod
\end{equation}
Under this identification, the standard inner product on $\bbC^4$ is
\begin{equation}
\DOT{x}{y} = \half\trace x\adj(y) \spacecomma\quad \adj
\begin{smatrix}a & b\\c & d\end{smatrix}
\coloneq
\begin{smatrix}d & -b\\-c & a\end{smatrix}
\end{equation}
and $\matSU{2}{}\subset\matSL{2}{\bbC}$ is identified with
$\bbR^4\subset\bbC^4$. In particular, for $X,\,Y\in \matsl{2}{\bbC}$
it holds $\DOT{x}{y}=-\half\trace xy$.

In order to treat irreducibility, the following lemma translates the
notion of eigenline to a more convenient form.  With
$\begin{smatrix}a\\b\end{smatrix}^\perp \coloneq [\begin{smallmatrix}-b &
      a\end{smallmatrix}]$ consider the double cover
\begin{equation}
  \label{eq:double-cover-u}
  \hat x\in\bbC^2\setminus\{0\}\to\{x\in\matsl{2}{\bbC}\suchthat \det
  x = 0\} \spacecomma\quad \hat x\mapsto x= \hat x \hat x^\perp
  \spaceperiod
\end{equation}

\begin{lemma}
  \label{lem:eigenvalue}
  \theoremname{Eigenvalues} $\ell\in\bbC^2\setminus\{0\}$ is an
  eigenvector of the invertible matrix $x\in\matSL{2}{\bbC}$ if and
  only if $\DOT{x}{\ell\ell^\perp} = 0$.
\end{lemma}

\begin{proof}
  For any $p,\,q\in\bbC^2\setminus\{0\}$,
  \begin{equation}
  \trace q p^\perp = p^\perp q = \det(p,\,q) \spaceperiod
  \end{equation}
  So with $p = \ell$, $q = x\ell$, and $y=\ell\ell^\perp$
  \begin{equation}
  2\DOT{x}{y} = 2\DOT{x}{\ell\ell^\perp} = \det(x\ell,\,\ell)
  \end{equation}
  so $\DOT{x}{y}=0$ if and only if $x\ell$ and $\ell$ are dependent,
  that is, if and only if $\ell$ is an eigenline of $x$.
\end{proof}

Let $\calP$ be the group generated by
$P_0=\id,\,P_1,\,P_2,\,P_3\in\matSL{2}{\bbC}$.  Under the above
identification $\bbC^4\cong\matgl{2}{\bbC}$ let
$P=(P_0,\,P_1,\,P_2,\,P_3)\in\matN{4}{\bbC}$ be a matrix with columns
$P_k\in\bbC^4$.  Let $T = \transpose{P}P\in\matsym{4}{\bbC}$, so
$T_{ij} = \DOT{P_i}{P_j}$.

\begin{theorem}
  \label{thm:unitarizability}
  \theoremname{Irreducibility and unitarizability} With $P$ and $T$ as
  above,
  \begin{enumerate}
  \item
    \label{item:unitarizability-1}
    $\calP$ is irreducible if and only if $\rank T\ge 3$.
  \item
    \label{item:unitarizability-2}
    Assuming \eqref{item:unitarizability-1}, $\calP$ is $\matSU{2}{}$
    unitarizable if and only if $T$ is real positive semidefinite.
  \end{enumerate}
\end{theorem}

\begin{proof}
  We have the factorization
  \begin{equation}
    T =
    \begin{bmatrix}\id & 0\\ \transpose{V} & \id\end{bmatrix}
    \begin{bmatrix}\id & 0\\ 0 & \transpose{X}\end{bmatrix}
    \begin{bmatrix}\id & 0\\ 0 & X\end{bmatrix}
    \begin{bmatrix}\id & V\\ 0 & \id\end{bmatrix}
   \spacecomma\quad P = \begin{bmatrix}1 & V\\0 & X\end{bmatrix}
     \spaceperiod
  \end{equation}
  To prove~\eqref{item:unitarizability-1}, let
  $X=(x_0,\,x_1,\,x_2)\in\matN{3}{\bbC}$ be the lower right $3\times
  3$ submatrix of $P$, that is the matrix with columns given by the
  tracefree parts of $P_1,\,P_2,\,P_3$, and let $Y\coloneq
  \transpose{X}X$.  By \cref{lem:eigenvalue}, $\calP$ is irreducible
  if and only if $\ker \transpose{X}\cap\calL = \{0\}$.  By
  \cref{lem:c3geometry} this is if and only if $\rank Y \ge 2$.  Since
  $\rank T = 1 + \rank Y$, this is if and only if $\rank T\ge 3$.

  To prove~\eqref{item:unitarizability-2}, if $\calP$ is $\matSU{2}{}$
  unitarizable, it may be assumed without loss of generality that
  $\calP\subset \matSU{2}{}$.  Then $P\in\matN{4}{\bbR}$, so $T =
  \transpose{T}T\in\matsym{4}{\bbR}$ is real positive semidefinite.

  Conversely, if $T$ is real positive semidefinite, then $W =
  \transpose{X}X = \transpose{Y}Y$ for some
  $Y\in\matN{3}{\bbR}$. Replacing $Y\mapsto \NEG Y$ if necessary, then
  $\det X = \det Y$, so by
  \cref{lem:c3geometry}\eqref{item:c3geometry-2}, there exists
  $S\in\matSO{3}{\bbC}$ such that $X = S Y$.  Let
  $C\in\matSL{2}{\bbC}$ be a lift of $S$ via the double cover
  $\matSL{2}{\bbC}\mapsto \matSO{3}{\bbC}$ defined with respect to
  $\id,\,e_1,\,e_2,\,e_3$.  Note that this double cover is given by
  conjugation on $\matsl{2}{\bbC}\cong \bbC^3$.  Then $C$ unitarizes
  $\calP$.
\end{proof}

\subsection{Constructing the surface}

Once the potential for a surface is obtained via the unitary flow, the
surface is constructed as follows:
\begin{itemize}
\item
  Compute the unitarizer of the monodromy.
\item
  Compute curvature lines.
\item
  Compute the fundamental piece of the surface via the DPW
  construction
\item
  Build the surface from the fundamental piece by reflections.
\end{itemize}

\subsubsection{The unitarizer}
Due to the symmetry \eqref{symmetric-potential} of the potential and
of the monodromies with basepoint $z=1$ the unitarizer is diagonal and
can be computed as follows.  With the notation $a^\ast(\lambda) =
\ol{a(1/\ol{\lambda})}$ write
\begin{equation}
M_0 = \begin{bmatrix}a & b \\ c & a^\ast\end{bmatrix} \spaceperiod
\end{equation}
By the unitarizability of $M_0$ by a diagonal loop, $p
\coloneq-c^\ast/b = -c/b^\ast$ takes values in $\bbR_+$ along $\bbS^1$
away from its zeros and poles, which are even.  Let $f =
\prod(\lambda-\alpha)/\prod(\lambda-\beta)$ so that $f^\ast f$ has the
same zeros and poles as $p$.  Then $q = p/(f^\ast f)$ takes values in
$\bbR_+$ along $\bbS^1$ without zeros or poles.  Let $y^{\ast}y = q$
be the scalar $\matGL{1}{\bbC}$ Birkhoff factorization, so $y$ is
holomorphic in the unit disk.  Then with $x = fy$ the loop $\diag(
{x}^{1/2},\, {x}^{-1/2} )$ is the required unitarizer, holomorphic on
the open unit disk.

\subsubsection{Curvature lines}
Let $Q=q(z)\deriv z^2$ the Hopf differential of the CMC surface.  The
curvature line coordinate $v$ satisfies $\deriv v^2 = Q(z) \deriv
z^2$.  Curvature line coordinates can be computed by computing
$\int\!\!\sqrt{Q(z)}\,\deriv z$ over the domain.

The surface is computed numerically by dividing the domain into
polygons (triangles or quadrilaterals) and mapping via the CMC
immersion these triangles to $\bbR^3$.  In the computation of
curvature lines described above, the polygon edges are unrelated to
the curvature lines.

Quadrilaterals whose edges are along curvature lines can be computed
as follows.  Divide the domain into quadrilaterals whose edges are
curvature lines and such that the umbilics are at corners of the
quadrilaterals.  For each quadrilateral, pull back the potential to
curvature line coordinates.

This computation is complicated by the fact that the maps from
curvature lines rectangles to the domain are singular at the umbilics,
and the potential is singular at the umbilics.  The potential can be
desingularized locally at an umbilic $z_0$ by a coordinate change of
the form $z = z_0 + w^n$ and a gauge.

\vspace{0.0975cm}\begin{textblock*}{16cm}(14.5cm,1cm)
\begin{minipage}{5cm}
\begin{figure}
  \ximage[1]
    {s3-platonic-altoct.pdf}
    {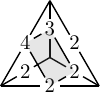}
    {\scriptsize (\textsc{a}) alternate octahedron \vspace{0.35cm}}
  \ximage[1]
    {s3-platonic-ico.pdf}
    {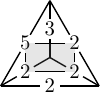}
    {\scriptsize (\textsc{b}) icosahedron \vspace{0.35cm}}
  \ximage[1]
    {s3-platonic-dodec.pdf}
    {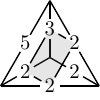}
    {\scriptsize (\textsc{c}) dodecahedron \vspace{0.35cm}}
  \ximage[1]
    {s3-platonic-altico.pdf}
    {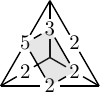}
    {\scriptsize (\textsc{d}) alternate icosahedron \vspace{0.35cm}}
  \captionsetup{margin=0cm}
  \caption{\label{fig:s3-platonic-b}Minimal surfaces in $\bbS^3$ with Platonic symmetry}
\end{figure}
\end{minipage}
\end{textblock*}

\subsubsection{Building the surface}
In general the position of the surface in space is not controlled, so
to build the surface it must first be put into a standard position,
where a group of standard reflections can be applied.  To do so,
compute the four generating reflections $R_k$ in the isometry group
$\Iso\,\bbR^3$ of $\bbR^3$.  Conjugate them to standard reflections
$S_k$ via $CR_kC^{-1} = S_k$.  Then the surface after being moved via
$x\mapsto Cx$ has the standard reflections $S_k$ as symmetries.

\subsubsection{The bulge count for families of CMC surfaces}
The surfaces constructed in this paper allow for non-trivial 1-parameter deformations 
within the space of CMC surfaces with the same combinatorics.
A natural question, also considered by \cite{Grosse-Brauckmann_1993}, is whether different surfaces  with the same combinatorics,
but  which swap neck and bulge,
belong to the same family of CMC surfaces, for example \cref{fig:r3-lattice3-a}\textsc{a-b}, and \cref{fig:r3-lattice2-a}\textsc{a-b}.
It turns out that these  examples belong to different families.
We denote by a leg of the surface a cylindrical piece obtained from 
the trajectories of the Hopf differential, that is from the curvature line parametrisation.
Although the images are labeled according to whether there are bulges or necks where the legs meet,
in this section we rather count the number of bulges on each leg.

We show that this number is an invariant
in the case of surfaces without ends.
Note that in this case, there is a covering $\Sigma\to\CPone$
by a compact Riemann surface $\Sigma$ on which the pullback of the DPW potential has only apparent
singularities. Phrased differently, $\Sigma$ is the surface on which the first and second fundamental forms 
are well-defined and smooth, that is for compact CMC surfaces in the 3-sphere, $\Sigma$ is just the underlying Riemann surface, and
in the case of periodic CMC surfaces in $\bbR^3$, $\Sigma$ is the  Riemann surface quotient of the CMC surface by the translational symmetries.

We  construct surfaces starting from Delaunay cylinders by deforming
the eigenvalues $\nu_i.$ In the case of cylinders without umbilics, all four eigenvalues are $\nu_i=\tfrac{1}{4}.$
At the starting point, $\Sigma$ is a torus and the relevant moduli space of flat connections $\nabla$ on $\Sigma$ has only reducible points. The underlying holomorphic bundles
(equipped with the $(0,1)$-parts $\bar\del^\nabla$  of the connections) are semistable, i.e.~if  they admit holomorphic
line subbundles of degree 0. A holomorphic structure (on a rank 2 bundle over a compact Riemann surface of degree 0) is called unstable if there exist a holomorphic line subbundle of positive degree
and they are  called stable if every  holomorphic line subbundle has negative degree.
This notion is relevant to us  since an unstable holomorphic structure does not admit
a flat unitary connection. 
Spectral parameters $\lambda$ at which the holomorphic structure  is unstable are isolated
in the spectral plane.  Moreover, for CMC surfaces based on quadrilaterals, the number of those values of spectral parameters within 
a bounded region is always finite and can only change during a deformation by values crossing the boundary of that region.
Values of the spectral parameter at which the holomorphic structure is 
unstable  cannot cross
the unit circle, as the connections on the unit circle are unitary. 
For the initial torus the bundle is semistable for all spectral values.
The number of values at which the holomorphic bundle  becomes 
unstable within infinitesimal deformation of the eigenvalues $\nu_i$ can be identified with
the number of bulges on the leg of the initial 
Delaunay cylinder;  for more details  see \cite{Heller_Heller_Schmitt_2018}.
Actually this number coincides with the number of zeros of the holomorphic function $x$  in \eqref{symmetric-potential} inside the unit circle; see also
\cite{Heller_Heller_2016,Heller_Schmitt_2015}

\typeout{== section4 ============================================}\section{Fuchsian DPW potentials}\label{sec:FuchsianDPW4}

The aim of this section is to prove the existence of Fuchsian DPW
potentials of the form~\eqref{eq:potential} for CMC quadrilaterals
without Delaunay ends.  This generalizes previous work by the second
author~\cite{Heller_2013} for the Lawson genus 2 surface.  Similar
results have been obtained by Manca~\cite{Manca_2020}. Our arguments
are more geometric and prove the existence of a Fuchsian potential on
a 4-punctured sphere for all surfaces obtained by CMC quadrilaterals.

\vspace{0.0975cm}\begin{textblock*}{16cm}(14.5cm,7.5cm)
\begin{minipage}{5cm}
\begin{figure}
  \ximage[1]
    {s3-cell-5.pdf}
    {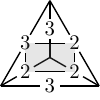}
    {\scriptsize (\textsc{a}) $5$-cell \vspace{0.35cm}}
  \ximage[1]
    {s3-cell-16.pdf}
    {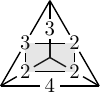}
    {\scriptsize (\textsc{b}) $16$- or $8$-cell \vspace{0.35cm}}
  \ximage[1]
    {s3-cell-600.pdf}
    {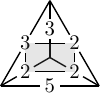}
    {\scriptsize (\textsc{c}) $600$- or $120$-cell \vspace{0.35cm}}
  \captionsetup{margin=0cm}
  \caption{\label{fig:s3-cell-a}Minimal surfaces in $\bbS^3$ with $n$-cell symmetry}
\end{figure}
\end{minipage}
\end{textblock*}

\subsection{Setup}
Let $f\colon \Sigma\to M$ (where $M\in\{\bbS^3,\bbR^3\}$) be a
complete CMC surface without Delaunay ends. Assume that $f$ is build
from a fundamental piece $P$ by the group $G$ generated by the
reflections across totally geodesic subspaces along geodesic arcs
contained in $P$. Assume that $P$ has the topology of a (closed) disc.

The surface $f$ is equivariant with respect to the (discrete) group
$G$ acting on $\Sigma$ by conformal transformations and on the ambient
space $M$ by a representation $\rho$ into the space of isometries. Let
$G^o\subset G$ be the subgroup of orientation preserving
(i.e. holomorphic) symmetries on $\Sigma$.

\subsection{Local theory}
The first step in our derivation of a Fuchsian DPW potential is the
converse of \cref{thm:cmc-polygon-vertex}. This means that at fixed
points of a rotational symmetry there always exists DPW potentials
with Fuchsian singularity on the quotient.

Let $p\in\Sigma$ be a fixed point of some rotation given by an element
in $G^o$. Then there exists $k\in \bbN$ and $g\in G^o$ of order $k$
such that $g(p)=p$ and such that for any $h\in G^o$ with $h(p)=p$
there exists $l\in \bbN$ with $g^l=h$.

\begin{lemma}
There exists $D\in \matSU{2}{}$ of order $2k$ and a local DPW
potential $\eta$ for $f$ on an open $g$-invariant neighbourhood of $p$
such that
\begin{equation}
g^*\eta=D\eta D^{-1}\spaceperiod
\end{equation}
\end{lemma}

\begin{proof}
  Consider Dorfmeister's normalized potential (see for
  example~\cite{Wu_1999}) which takes the form
\begin{equation}
  \eta^{nor}=\begin{bmatrix}0& \lambda^{-1}
  f(z,0)\\ \frac{q}{f(z,0)}&0\end{bmatrix}dz
\end{equation}
where $z$ is a local holomorphic coordinate centered in $p$ such that
$g^*z=e^{\frac{2\pi i}{k}}z$, $Q=q(dz)^2$ is the Hopf differential and
$f(z,w)$ is a holomorphic function such that $f(z,\bar z)dz d\bar z$
is the induced metric of the surface. As $g^*dz=e^{\frac{2\pi
    \imi}{k}}dz$ and $g^*d\bar z=e^{\frac{2\pi i}{k}}d\bar z$ the
result follows.
\end{proof}

\begin{proposition}
\label{potlocregsing}
There exists  a local meromorphic DPW potential of $f$ on
$\Sigma/G^o$ with a Fuchsian singularity at $p\mod G^o$.  The
eigenvalues of the residue are $\pm\frac{1}{2k},$ independently of
$\lambda$, where $k$ is the order of the stabilizer group of $p$.

Likewise, there  exists a local meromorphic DPW potential of $f$
on $\Sigma/G^o$ with Fuchsian singularity at $p\mod G^o$ such that the
eigenvalues of the residue are $\pm\frac{ k-1}{2k}$.\end{proposition}
\begin{proof}
Consider $w=z^k$ which is a holomorphic coordinate centered at $p\mod
G^o\in \Sigma/G^o$. Consider the positive gauge $e=\diag{(\sqrt{z}
  ,\frac{1}{\sqrt{z}})}$ of $\spin$ $-1$. Then
\begin{equation}
(d+\eta^{nor}).e=d+ e^{-1} de +e^{-1} \eta^{nor}e
\end{equation}
is a well-defined meromorphic DPW potential with apparent Fuchsian
singularity at $p$. As this potential is clearly invariant under
pull-back by $g$ we have proven the first part of the proposition.

For the second part and $k=2l+1$ consider the gauge $\tilde
e=\diag{(z^{-l},z^l)}$ while for $k=2l$ consider the gauge $\tilde
e=\diag{(z^{-l+\tfrac{1}{2}},z^{l-\tfrac{1}{2}})}$, and proceed as in
the first part of the proof.
\end{proof}

\subsection{Global theory}

Our aim is to construct a DPW potential on the $\Sigma/G^o$.  Recall
that by assumption the fundamental piece $P$ of the Riemann surface
$\Sigma$ is of the topological type of a disc.

\begin{lemma}
The Riemann surface $\Sigma/G^o$ is the projective line.
\end{lemma}
\begin{proof}
By the Riemann mapping theorem there exists a holomorphic map from $P$
to the unit disc. Schwarzian reflection yields a holomorphic map from
$\Sigma$ to $\CPone$, branched at the fixed points of $G^o$. By its
construction, this map is invariant under $G^o$.
\end{proof}
For simplicity of the arguments, we will assume that $n$ is even in
the following.

\begin{lemma}\label{lem:uniP1pot}
Let $n$ be even.  There  exists a unitary potential $\mu$ on the
$n$-punctured Riemann sphere such that
\begin{itemize}
\item $\mu$ is singular exactly at the branch values of
  $\Sigma\to\Sigma/ G^o=\CPone;$
\item the pull-back of $\mu$ generates $f$ on the covering $\Sigma$.
\end{itemize}
\end{lemma}
\begin{proof}
Let $\{z_1,\dots,z_n\}\subset \CPone$ be the branch values of
$\Sigma\to\Sigma/ G^o=\CPone$ and $S\subset\Sigma$ its preimage.
Denote the reflection planes of the fundamental piece by $P_0,\dots,
P_{n-1},$ with outward oriented unit normals $N_0,\dots , N_{n-1},$
respectively, such that
\begin{equation}
z_m\in P_{m-1}\cap P_m\quad \forall m\in\{1,\dots,n\}\spaceperiod
\end{equation}
Denote by $g_m$ the compositions of the reflection across $P_{m-1}$
and $P_m$. Then $G^0$ is generated by
$\{g_m\mid\,m\in\{1,\dots,n\}\}$.

Let $M$ be euclidean 3-space or the 3-sphere, and let $d=3$ and $d=4$
accordingly, so that $\Iso(M)=\matSO{4}{\bbR}$ or
$\Iso(M)=\matSO{3}{\bbR}\ltimes \bbR^3$.

Consider the group
\begin{equation}
H\subset \matSpin{n}\times \Iso(M)
\end{equation}
generated by the elements
\begin{equation}
 \hat g_m:=(N_m\cdot N_{m-1},g_m), \quad m=1,\dots,n
\end{equation}
where $\cdot$ denotes Clifford multiplication. This gives a group
extension
\begin{equation}
\{\text{id}\}\to\bbZ_2\to H\to G^0\to \{\text{id}\}\spaceperiod
\end{equation}
Note that $\hat g_m$ has order $2k$ if $g_m$ has order $k$. Similarly,
since $n$ is even, the product $\hat g_n\dots\hat g_1$ is trivial.
Consequently, we have a representation
\begin{equation}
h\colon\pi_1(\CPone\setminus\{z_1,\dots,z_n\},*)\to H\spaceperiod
\end{equation}

As
\begin{equation}
\Sigma\setminus S\to\CPone\setminus\{z_1,\dots,z_n\}
\end{equation}
  is a (unbranched) covering, the fundamental group of
  $\Sigma\setminus S$ (with appropriate base point) is a subgroup of
  the first fundamental group of $\CPone\setminus\{z_1,\dots,z_n\}$
  with corresponding base point. By construction, the induced
  representation of $\pi_1(\Sigma\setminus S,*)\to G^0$ is trivial,
  and the induced representation of $h$ takes values in
\begin{equation}
\bbZ_2=\bbZ_2\times\{\text{id}\}\subset \matSpin{d}\times \Iso(M)
\end{equation}
such that a simple closed curve around any one of the points in $S$ is
mapped to the non-trivial element in $\bbZ_2$.

\vspace{0.0975cm}\begin{textblock*}{16cm}(14.5cm,7.5cm)
\begin{minipage}{5cm}
\begin{figure}
  \ximage[1]
    {s3-cell-24.pdf}
    {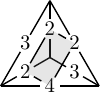}
    {\scriptsize (\textsc{a}) $24$-cell \vspace{0.35cm}}
  \ximage[1]
    {s3-cell-16-b.pdf}
    {diagram/tet-i-s3-cell-16.pdf}
    {\scriptsize (\textsc{b}) $16$- or $8$-cell \vspace{0.35cm}}
  \ximage[1]
    {s3-cell-600-b.pdf}
    {diagram/tet-i-s3-cell-600.pdf}
    {\scriptsize (\textsc{c}) $600$- or $120$-cell \vspace{0.35cm}}
  \captionsetup{margin=0cm}
  \caption{\label{fig:s3-cell-b}Minimal surfaces in $\bbS^3$ with $n$-cell symmetry}
\end{figure}
\end{minipage}
\end{textblock*}

By Riemann surface covering theory, we obtain a 2-fold covering
\begin{equation}
\hat\Sigma\to \Sigma\spacecomma
\end{equation}
branched over the points in $S$, with an action of $H$ by holomorphic
automorphisms on $\hat\Sigma$ such that
\begin{equation}
\hat\Sigma\to\hat\Sigma/H=\CPone
\end{equation}
is branched over $\{z_1,\dots,z_n\}$. Denote its preimage of $S$ by
$\hat S\subset\hat\Sigma$. Note that $H$ acts faithfully on
$\hat\Sigma\setminus \hat S$.

Consider the pull-back $\omega$ on $\hat\Sigma$ of the unitary
potential $\eta=F^{-1}\deriv F$ of $ f$. Note that, for minimal
$f\colon \Sigma\to\mathbb S^3$ the unitary potential is given by
\begin{equation}
\label{eq:unipot}\eta=\lambda^{-1}\Phi+\Phi-\Phi^*-\lambda\Phi^*
\end{equation}
where
\begin{equation}
\Phi=\tfrac{1}{2}(f^{-1}df)^{1,0}\quad\text{and}\quad
\Phi^*=\tfrac{1}{2}(f^{-1}df)^{0,1}
\end{equation}
and similarly for CMC surfaces $f\colon \Sigma\to\mathbb R^3$.  Let
$\pi\colon H\to \matSU{2}{}$ be the projection to the rotational part
of the symmetry. From the construction~\eqref{eq:unipot} of the
unitary potential,
\begin{equation}
h^*\eta_\lambda=\eta_\lambda.\pi(h)
\end{equation}
for a holomorphic automorphism $h\in H$ (where, on the right hand
side, the gauge action of the constant matrix $\pi(h)$ is given by
conjugation).

\vspace{0.0975cm}\begin{textblock*}{16cm}(14.5cm,20.2cm)
\begin{minipage}{5cm}
\begin{figure}
  \ximage[1]
    {s3-torus-8.pdf}
    {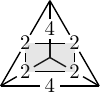}
    {\scriptsize (\textsc{a}) torus 8 \vspace{0.35cm}}
  \captionsetup{margin=0cm}
  \caption{\label{fig:s3-torus}Minimal torus in $\bbS^3$ with Delaunay ends}
\end{figure}
\end{minipage}
\end{textblock*}

Consider the free action of $H$ on $\hat \Sigma\setminus\hat
S\times\bbC^2$ given by
\begin{equation}
(p,v).h=(p.h,\pi(h^{-1}) (v))\spaceperiod
\end{equation}
The quotient
\begin{equation}
V=(\hat \Sigma\setminus\hat S\times\bbC^2)/H
\end{equation}
is a trivial smooth vector bundle of rank 2 over $\CPone\setminus
\{z_1,\dots,z_n\}$. We claim that the unitary potential $\omega$
yields a well-defined potential $\mu$ on this vector bundle: in fact,
the connection 1-form acts on $[p,v]\in(\hat \Sigma\setminus\hat
S\times\bbC^2)/H$ as
\begin{equation}
[p,\omega_p (v)]
\end{equation}
which is well-defined since
\begin{equation}
(p,\omega_p (v)).h=(ph,h^{-1} \omega_p h h^{-1}(v))=(ph, \omega_{ph}
  (h^{-1}v))\spaceperiod
\end{equation}
\end{proof}

\begin{proposition}
  Let $n$ be even.  There exist a meromorphic DPW potential $\xi$ on
  $\CPone$ with simple poles at $z_1,\dots,z_n$ and possible apparent
  singularity at $z=\infty$.
\end{proposition}

\begin{proof}
From \cref{lem:uniP1pot} we obtain a unitary potential $\mu$ on the
$n$-punctured sphere.  Let $l\in\{1,\dots,n\}$. By
\cref{potlocregsing} there exist a DPW gauge locally well-defined on a
punctured disc around $z_l$ which gauges $\mu$ into a meromorphic
potential with a Fuchsian singularity at $z_l$. Of course, the
holomorphic structures (i.e. the $(0,1)$-part) of a meromorphic
potential extends to the singular points.  Note that these gauges are
well-defined (i.e. have $\spin$ 1) as we have chosen the
representation $s\colon\pi_1(\Sigma\setminus S,\bbZ_2)$ to have local
monodromy $-1$ around every point in $S$.

Using these gauges as cocycles, we obtain a holomorphic
$\bbC^*$-family of flat $\matSL{2}{\bbC}$-connections $d+\hat\mu$ with
the following properties:
\begin{itemize}
\item the induced family of holomorphic structures extends to
  $\lambda=0$ to give a holomorphic rank 2 bundle $\mathcal
  E_0\to\CPone$ with trivial holomorphic determinant;
\item the connections $d+\hat\mu$ have Fuchsian singularities with
  $\lambda$-independent eigenvalues $\pm\tfrac{1}{2k};$
\item the complex linear part of the family of connections has a first
  order pole at $\lambda=0$, i.e., $\lambda\mapsto\lambda
  (\hat\mu)^{1,0}(\lambda)$ extends to $\lambda=0$.
\end{itemize}
Note that all bundles have trivial determinant. Hence, by the
Birkhoff-Grothendieck theorem, the bundle type of $\mathcal E_0$ is
$\mathcal O(d)\oplus\mathcal O(-d)$ for some $d\in\bbN$.

First consider the case $d=0$. Then, the bundle type $\mathcal
E_\lambda$ is locally constant on an open disc near $\lambda=0$. In
particular, there exists a smooth positive family of gauge
transformations $g_\lambda$ (holomorphic in $\lambda$) such that
\begin{equation}
((d+\hat\mu).g_\lambda)^{0,1}=d^{0,1}
\end{equation}
is the trivial holomorphic structure on the rank 2 vector bundle
$\bbC^2\to\CPone$. Thus
\begin{equation}
  d+\xi:=(d+\hat\mu).g_\lambda
\end{equation}
is the meromorphic DPW potential which has only Fuchsian singularities
and no apparent singularity at $\infty$.

Let $d>0$. Let $z\colon\CPone\setminus\{\infty\}\to\bbC$ be an affine
holomorphic coordinate and assume without loss of generality that
$z_l\neq\infty\quad\forall\; l\in\{1,\dots,n\}$.  There exists an
integer $0\leq s\leq d$ such that on a punctured disc $\mathcal
D\setminus\{0\}$ around $\lambda=0$ all bundles are of the holomorphic
type $\mathcal O(s)\oplus\mathcal O(-s)$.  By the family version of
Birkhoff-Grothendieck, there exists a holomorphic function
$r\colon\mathcal D\to\bbC$ with $r(0)=0$ such that the holomorphic
bundle $((d+\hat\mu).g_\lambda)^{0,1}$ has the cocycle (for the
covering $\mathcal U_+:=\bbC,\mathcal U_-:=\bbC
P^1\setminus\{\infty\}$ of $\CPone$)
\begin{equation}
\label{bundlenormalform}
\begin{bmatrix} z^{-d}& z^s r(\lambda)\\ 0& z^d\end{bmatrix}
  =
  \begin{bmatrix} r(\lambda) z^s&0\\z^d&\tfrac{1}{r(\lambda)}z^{-s}\end{bmatrix}
  \begin{bmatrix}\tfrac{1}{r(\lambda)}z^{-d-s}&1\\-1&0\end{bmatrix}\spacecomma
\end{equation}
where the equality obviously holds only for $r(\lambda)\neq0$. Again,
there exists a DPW gauge $g_\lambda$ which gauges $\mathcal E_\lambda$
into the above form~\eqref{bundlenormalform}.  This means that there
exists a pair $(g^+_\lambda,\,g^-_\lambda)$ of DPW gauges on $\mathcal
U_+$ respectively $\mathcal U_-$ which differ by the
gauge~\eqref{bundlenormalform} and gauge $(d+\hat\mu)^{0,1}$ on
$\mathcal U_{\pm}$ to the trivial holomorphic structure on
$\bbC^2\to\mathcal U_\pm$. Then,
\begin{equation}
(d+\hat\mu).g^+_\lambda=:d+\xi
\end{equation}
yields the meromorphic potential $\xi$ with Fuchsian singularities at
$z_k$ and an apparent singularity at $\infty$.
\end{proof}

\subsubsection{CMC quadrilaterals}
Finally, we consider the case of CMC quadrilaterals, i.e., $n=4$. We
show that these are always determined by a Fuchsian DPW
potential~\eqref{eq:potential}, which, assuming an additional
symmetry, is of the form~\eqref{symmetric-potential}.

\begin{lemma}
For $n=4$ the bundle type of $\mathcal E_0$ is either trivial or
$\mathcal O(1)\oplus\mathcal O(-1)\to\CPone$.
\end{lemma}

\begin{proof}
Assume the bundle type at $\lambda=0$ is $\mathcal O(d)\oplus\mathcal
O(-d)\to\CPone$ for some $d>1$. The Higgs field
$\Phi:=\text{res}_{\lambda=0}\xi$ is a meromorphic section of the
bundle
\begin{equation}
K\otimes \text{End}_0(\mathcal E_0)\to\CPone
\end{equation}
where $K=\mathcal O(-2)$ is the canonical bundle of $\CPone$ and
$\text{End}_0(\mathcal E_0)$ denotes the trace-free endomorphisms of
$\mathcal E_0$.  Moreover, $\Phi$ is nilpotent as the immersion is
conformal, has at most simple poles at $z_1,\dots,z_4$ by construction
and does not vanish on $\CPone$ as $f$ is an immersion. Using the
decomposition $\mathcal E_0=\mathcal O(d)\oplus\mathcal O(-d)$ the
Higgs field is of the form
\begin{equation}
\Phi=\begin{bmatrix} a& b\\ c&-a\end{bmatrix}
\end{equation}
where $a,b,c$ are meromorphic sections in $\mathcal O(-2)$, $\mathcal
O(-2+2d),\mathcal O(-2-2d),$ respectively, with at most simple poles
at $z_1,\dots,z_4$. Hence $c=0$. As $\Phi$ is nilpotent $a=0$ as
well. For $d>1$, $-2+2d>0$ and $b$ would have a zero contradicting the
fact that $\Phi$ is nowhere vanishing.
\end{proof}

\begin{theorem}
  \label{fuchsdpw}
Let $f$ be a complete CMC surface without Delaunay ends in $\bbS^3$ or
$\bbR^3$.  If $f$ is built from a CMC quadrilateral in a fundamental
tetrahedron of a tessellation of the ambient space then it is obtained
from a Fuchsian DPW potential~\eqref{eq:potential} on the 4-punctured
sphere.
\end{theorem}

\begin{proof}
We give a proof by contradiction.  Assume that the bundle type at
$\lambda=0$ is
\begin{equation}
\label{bt1m1}\mathcal E_0=\mathcal O(-1)\oplus\mathcal
  O(1)\spaceperiod
\end{equation}
  By the proof of \cref{potlocregsing} the nilpotent
  $\lambda^{-1}$-part $\Phi=\xi_{-1}$ of the meromorphic potential has
  no zeros, and poles of order $1$ at the 4 branch points
  $z_1,\dots,z_4$. Thus, with respect to \eqref{bt1m1}, it must be of
  the form
\begin{equation}
\Phi=\begin{bmatrix} 0 & s_{-D}\\0&0\end{bmatrix}
\end{equation}
where $s_{-D}\in\mathcal M(\CPone,\mathcal O(-4))$ is the unique
meromorphic section (up to scaling) with simple poles at $D=z_1+\dots
+z_4$. Moreover, the positive eigenvalues $\nu_i$ of the residues of
the connections are contained in the respective kernels of the
residues of $\xi_{-1}$. This equips $\mathcal E_0$ with a parabolic
structure (see for
example~\cite{Mehta_Seshadri_1980,Simpson_1990,Pirola_2007,Heller_Heller_2016,Biswas_Dumitrescu_Heller_2020}
for definitions and further references) which is unstable.  We denote
the parabolic bundle also by $\mathcal E_0$.  The pair $(\mathcal
E_0,\Phi)$ is a stable strongly parabolic Higgs pair.  Note that
\begin{equation}
\sum_i\nu_i<1\spaceperiod
\end{equation}
It is easy to see (compare with~\cite{Heller_Heller_2016}) that
$(\mathcal E_0,\Phi)$ is the only stable strongly parabolic Higgs pair
with nilpotent Higgs field on the 4-punctured sphere with unstable
underlying parabolic bundle.  Consider the compact Riemann surface
$X\to\CPone$ on which the rotational symmetry is trivial. Its Fuchsian
monodromy (given by uniformization) corresponds by the
Hitchin-Kobayashi correspondence to a stable nilpotent Higgs pair
\begin{equation}
  \biggl(S^*\oplus
  S\spacecomma\ \begin{bmatrix}0&1\\0&0\end{bmatrix}\biggr)\spacecomma
\end{equation}
where $S^2=K_X$. Its underlying holomorphic structure is unstable. As
the rotational symmetries act on $X$ we obtain, in the same manner as
for $f$, an strongly parabolic nilpotent Higgs pair with underlying
parabolic structure.  As the holomorphic structure is unstable, the
parabolic structure must be unstable as well; see~\cite{Biswas_1997},
and hence it must be $(\mathcal E_0,\Phi)$. Thus, the holomorphic
Higgs pair of $f$, i.e., $(\partial^\nabla,\Phi)$ would be gauge
equivalent to $(S^*\oplus
S,\bigl[\begin{smallmatrix}0&1\\0&0\end{smallmatrix}\bigr])$. This is
only possible if the Hopf differential of the minimal (respectively
CMC) surface vanishes (compare with~\cite[sections 2 and
  3]{Heller_2013}), which gives a contradiction.
\end{proof}

Finally, we show under which conditions the Fuchsian DPW potential
$\xi$ can be gauged into the form~\eqref{symmetric-potential}.

Note that a Fuchsian potential for a CMC quadrilateral defining a
compact embedded CMC surface cannot be adapted if all the 4 positive
eigenvalues of the residues are contained in $(0,\tfrac{1}{4})$.

\begin{corollary}\label{fuchsdpw2}
Assume that the potential of \cref{fuchsdpw} has equal pairs of
eigenvalues.  Then, there exist a coordinate change and a gauge such
that the potential is of the form~\eqref{symmetric-potential}.
\end{corollary}

\begin{proof}
We only sketch the proof. Assume that the eigenvalues at $z_0$ and
$z_1$, respectively $z_2$ and $z_3$ are equal.  First, apply a
so-called flip gauge which flips the eigenvalues at $z_2$ and $z_3$ by
adding $\mp\tfrac{1}{2}$. This can be achieved by conjugating the
potential by a DPW gauge which is constant in $z$ such that the
residues at $z_2$ and $z_3$ are lower respectively upper triangular,
and then gauge with
$\diag(\sqrt{\tfrac{z-z_2}{z-z_3}},\sqrt{\tfrac{z-z_3}{z-z_2}})$.
Denote the residues of the potential $\tilde \xi$ obtained in this way
by $R_k,$ and find $T$ such that $R_3=
TD^{-1}T^{-1}R_2TDT^{-1}$. Then, $T^{-1}\tilde\xi T$ turns out to be
of the form~\eqref{symmetric-potential}.
\end{proof}

\bibliographystyle{amsplain}
\bibliography{references}

\end{document}